\pdfoutput=1
\RequirePackage{silence}
\documentclass[a4paper,svgnames]{amsart}
\usepackage[hmarginratio={1:1},vmarginratio={1:1},lmargin=60.0pt,tmargin=60.0pt]{geometry}

\usepackage{booktabs}
\allowdisplaybreaks

\usepackage[T1]{fontenc}
\usepackage{latexsym,exscale,amsfonts,amssymb,mathtools}
\usepackage{amsmath,amsthm,amsfonts,amssymb,amscd,textcomp,bbm}
\usepackage{mathrsfs,stackrel}
\usepackage{etoolbox}
\usepackage{tikz,tikz-cd}
\usetikzlibrary{matrix,arrows.meta,decorations.markings,shapes.multipart,decorations.pathreplacing,decorations.pathmorphing}
\usepackage{aliascnt}
\usepackage{ytableau}
\usepackage{xparse}
\usepackage{cite}

\usepackage{dynkin-diagrams}

\usepackage{array}
\newcolumntype{C}{>{$}c<{$}}

\usepackage{xcolor,colortbl}

\definecolor{mygray}{gray}{0.6}
\definecolor{mygraydark}{gray}{0.4}
\definecolor{mygraylight}{gray}{0.85}
\definecolor{spinach}{RGB}{46,139,87}
\definecolor{tomato}{RGB}{255,99,71}
\definecolor{orchid}{RGB}{143,40,194}
\definecolor{neon}{RGB}{77,77,255}
\definecolor{pumpkin}{RGB}{224,180,80}
\definecolor{citron}{RGB}{190,180,90}

\definecolor{lava}{RGB}{207,16,32}
\definecolor{cream}{RGB}{255,253,208}
\definecolor{verdigris}{RGB}{67,179,174}
\definecolor{Black}{RGB}{0,0,0}
\definecolor{mydarkblue}{RGB}{10,10,170}
\definecolor{darkspinach}{RGB}{20,70,20}
\definecolor{darktomato}{RGB}{155,40,30}
\definecolor{darkorchid}{RGB}{50,10,100}
\definecolor{darklava}{RGB}{150,8,16}

\usepackage{todonotes}

\usepackage{enumitem}
\setlist[enumerate]{itemsep=0.15cm,label=\emph{\upshape(\alph*)}}
\setlist[enumerate,2]{itemsep=0.15cm,label=\emph{\upshape(\roman*)}}
\setlist[enumerate,3]{itemsep=0.15cm,label=\emph{\upshape(\Alph*)}}

\let\emph\relax
\DeclareTextFontCommand{\emph}{\em}

\renewcommand{\dots}{\text{...}}
\renewcommand{\vdots}{\rotatebox{90}{\text{...}}}
\renewcommand{\ddots}{\raisebox{0.175cm}{\rotatebox{-45}{\text{...}}}}
\newcommand{\rddots}{\raisebox{-0.075cm}{\rotatebox{45}{\text{...}}}}
\newcommand{\placeholder}{{}_{-}}

\newcommand{\ie}{\text{i.e.}}
\newcommand{\eg}{\text{e.g.}}
\newcommand{\cf}{\text{cf.}}
\newcommand{\etc}{\text{etc.}}

\newcommand{\C}{\mathbb{C}}
\newcommand{\R}{\mathbb{R}}
\newcommand{\Rplus}{\mathbb{R}_{\geq 0}}
\newcommand{\N}{\mathbb{Z}_{\geq 0}}

\newcommand{\Z}{\mathbb{Z}}
\newcommand{\Zplus}{\mathbb{Z}_{\geq 0}}
\newcommand{\K}{\mathbb{K}}
\newcommand{\F}[1][p]{\mathbb{F}_{#1}}
\newcommand{\bF}[1][p]{\bar{\mathbb{F}}_{#1}}

\newcommand{\munit}{\mathbbm{1}}

\newcommand{\catstuff}[1]{\mathbf{#1}}

\newcommand{\obstuff}[1]{\mathtt{#1}}

\newcommand{\pfdimplain}{\mathrm{PFdim}}
\newcommand{\pfdim}{\mathrm{PFdim}^{f}}
\newcommand{\pfdimk}[1][k]{\mathrm{PFdim}_{#1}}
\newcommand{\pfdimtwo}{\mathrm{PFdim}^{VJ}}
\newcommand{\pfdimthree}{\mathrm{PFdim}^{f}}

\newcommand{\asymbol}{a}
\newcommand{\bsymbol}{b}

\newcommand{\lsymbol}{l}

\newcommand{\rep}{\catstuff{Rep}}

\newcommand{\sbim}{\catstuff{SBim}}

\newcommand{\simple}[1][1]{\mathrm{L}_{#1}}
\newcommand{\inde}[1][1]{\mathrm{Z}_{#1}}
\newcommand{\prin}[1][1]{\mathrm{P}_{#1}}

\usepackage[all]{xy}
\usepackage{tikz}
\usetikzlibrary{cd}
\usetikzlibrary{decorations}
\usetikzlibrary{decorations.markings}
\usetikzlibrary{decorations.pathreplacing}
\usetikzlibrary{decorations.pathmorphing}
\usetikzlibrary{arrows.meta,shapes,positioning,matrix,calc}
\usetikzlibrary{shapes.callouts}
\tikzset{anchorbase/.style={baseline={([yshift=-0.5ex]current bounding box.center)}},
tinynodes/.style={font=\tiny,text height=0.25ex,text depth=0.05ex},
smallnodes/.style={font=\scriptsize,text height=0.75ex,text depth=0.15ex},
}

\def\NewTheorem#1{%
\newaliascnt{#1}{equation}%
\newtheorem{#1}[#1]{#1}%
\aliascntresetthe{#1}%
\expandafter\def\csname #1autorefname\endcsname{#1}%
}
\def\equationautorefname~#1\null{(#1)\null}

\numberwithin{equation}{subsection}

\NewTheorem{Proposition}
\NewTheorem{Theorem}
\NewTheorem{Corollary}
\AtEndEnvironment{Corollary}{\null\hfill$\square$}%
\NewTheorem{Lemma}
\theoremstyle{definition}
\NewTheorem{Definition}
\AtEndEnvironment{Definition}{\null\hfill$\Diamond$}%
\NewTheorem{Classification Problem}
\AtEndEnvironment{Classification Problem}{\null\hfill$\Diamond$}%
\NewTheorem{Notation}
\AtEndEnvironment{Notation}{\null\hfill$\Diamond$}%
\NewTheorem{Convention}
\AtEndEnvironment{Convention}{\null\hfill$\Diamond$}%
\NewTheorem{Example}
\AtEndEnvironment{Example}{\null\hfill$\Diamond$}%
\NewTheorem{Examples}
\AtEndEnvironment{Examples}{\vskip-10mm\null\hfill$\Diamond$}%

\theoremstyle{remark}
\NewTheorem{Remark}
\AtEndEnvironment{Remark}{\null\hfill$\Diamond$}%
\NewTheorem{Question}
\AtEndEnvironment{Question}{\null\hfill$\Diamond$}%
\NewTheorem{Assumption}

\usepackage[hypertexnames=false]{hyperref}
\usepackage{bookmark}
\hypersetup{
pdftoolbar=true,
pdfmenubar=true,
pdffitwindow=false,
pdfstartview={FitH},
pdftitle={Asymptotics in infinite monoidal categories},
pdfauthor={Abel Lacabanne, Daniel Tubbenhauer and Pedro Vaz},
pdfsubject={},
pdfcreator={Abel Lacabanne, Daniel Tubbenhauer and Pedro Vaz},
pdfproducer={Abel Lacabanne, Daniel Tubbenhauer and Pedro Vaz},
pdfkeywords={},
pdfnewwindow=true,
colorlinks=true,
linkcolor=mydarkblue,
citecolor=teal,
filecolor=magenta,
urlcolor=orchid,
linkbordercolor=lava,
citebordercolor=teal,
urlbordercolor=orchid,
linktocpage=true
}

\def\makeautorefname#1#2{\csdef{#1autorefname}{#2}}
\makeautorefname{section}{Section}%
\makeautorefname{subsection}{Section}%
\makeautorefname{subsubsection}{Section}%

\begin{document}
\title[Asymptotics in infinite monoidal categories]{Asymptotics in infinite monoidal categories}
\author[A. Lacabanne, D. Tubbenhauer and P. Vaz]{Abel Lacabanne, Daniel Tubbenhauer and Pedro Vaz}

\address{A.L.: Laboratoire de Math{\'e}matiques Blaise Pascal (UMR 6620), Universit{\'e} Clermont Auvergne, Campus Universitaire des C{\'e}zeaux, 3 place Vasarely, 63178 Aubi{\`e}re Cedex, France,\newline \href{http://www.normalesup.org/~lacabanne}{www.normalesup.org/$\sim$lacabanne},
\href{https://orcid.org/0000-0001-8691-3270}{ORCID 0000-0001-8691-3270}}
\email{abel.lacabanne@uca.fr}

\address{D.T.: The University of Sydney, School of Mathematics and Statistics F07, Office Carslaw 827, NSW 2006, Australia, \href{http://www.dtubbenhauer.com}{www.dtubbenhauer.com}, \href{https://orcid.org/0000-0001-7265-5047}{ORCID 0000-0001-7265-5047}}
\email{daniel.tubbenhauer@sydney.edu.au}

\address{P.V.: Institut de Recherche en Math{\'e}matique et Physique, 
Universit{\'e} catholique de Louvain, Chemin du Cyclotron 2,  
1348 Louvain-la-Neuve, Belgium, \href{https://perso.uclouvain.be/pedro.vaz}{https://perso.uclouvain.be/pedro.vaz}, \href{https://orcid.org/0000-0001-9422-4707}{ORCID 0000-0001-9422-4707}}
\email{pedro.vaz@uclouvain.be}

\begin{abstract}
We discuss formulas for 
the asymptotic growth rate of the number of summands in tensor
powers in certain (finite or infinite) monoidal categories. 
Our focus is on monoidal categories with infinitely many indecomposable objects, with 
our main tools being generalized Perron--Frobenius theory alongside techniques from random walks.
\end{abstract}


\addtocontents{toc}{\protect\setcounter{tocdepth}{1}}

\maketitle

\tableofcontents


\section{Introduction}\label{S:Intro}


\subsection{Growth problems}

Let $\catstuff{C}$ be an
additive Krull--Schmidt monoidal category.
Let $\obstuff{X}\in\catstuff{C}$ be an object of $\catstuff{C}$.
We define 
\begin{gather*}
\bsymbol_{n}=\bsymbol_{n}^{\catstuff{C},\obstuff{X}}:=\#\text{indecomposable summands in $\obstuff{X}^{\otimes n}$ counted with multiplicities}.
\end{gather*}
Let $\bsymbol\colon\N\to\N,\bsymbol(n)=\bsymbol_{n}$ denote the associated function.

\begin{Notation}
We use $\bsymbol_{n}$ to denote the entries of the sequence of numbers 
$(\bsymbol_{n})_{n\in\N}$, $\bsymbol$ for the associated function
and $\bsymbol(n)$ for the evaluation of $b$.
We use a similar notation for other sequences and functions.
\end{Notation}

\begin{Question}\label{Q:IntroductionMainQuestion}
The main point of this paper is to address, and partially answer, the following questions:
\begin{enumerate}

\item What is the \emph{dominating growth} of $\bsymbol_{n}$?
\label{Eq:IntroductionMainQuestionA}

\item Can we get an \emph{asymptotic 
formula} $\asymbol\colon\N\to\N$ expressing $\bsymbol$, e.g.
\begin{gather*}
\bsymbol(n)\sim\asymbol(n),
\end{gather*}
where $\asymbol$ is ``nice'', $\sim$ denotes asymptotically equal?\label{Eq:IntroductionMainQuestionB}

\item Say we have found $\asymbol$ as in the previous point.
Can we bound the \emph{variance} or 
\emph{mean absolute difference} $|\bsymbol_{n}-\asymbol_{n}|$, or alternatively the convergence 
rate of $\lim_{n\to\infty}\bsymbol(n)/\asymbol(n)=1$?
\label{Eq:IntroductionMainQuestionC}

\end{enumerate}
We refer to \ref{Eq:IntroductionMainQuestionA}, \ref{Eq:IntroductionMainQuestionB}, and \ref{Eq:IntroductionMainQuestionC} as the \emph{growth problems} associated with $(\catstuff{C},\obstuff{X})$.
We also explore several statements along the same lines, which we will also call growth problems.
\end{Question}

\begin{Remark}\label{R:IntroductionMainQuestion}
The ordering in \autoref{Q:IntroductionMainQuestion} roughly reflects increasing difficulty. To draw an analogy:
if $b$ would be the prime counting function, then 
\ref{Eq:IntroductionMainQuestionA} could say that it grows essentially linearly, \ref{Eq:IntroductionMainQuestionB} could be the prime number theorem, and \ref{Eq:IntroductionMainQuestionC} 
could be (a consequence of) the Riemann hypothesis.
\end{Remark}

\begin{Example}
If $\obstuff{X}$ and the monoidal unit $\munit$ are indecomposable, and we have $\obstuff{X}\otimes\obstuff{X}\cong\munit\oplus\obstuff{X}$, then $(\bsymbol_{n})_{n\in\N}$ is the Fibonacci sequence (with the first two terms equal to $1$). In this case the dominating growth is $\phi^{n}$ 
for $\phi\approx1.618$ the golden ratio.
We have $b(n)\sim\frac{\phi}{\sqrt{5}}\cdot\phi^{n}$ and 
$|b_{n}-a_{n}|\leq(\phi-1)^{n}$ for $\phi-1\approx0.618$.

Note the following: the action matrix of $\placeholder\otimes\obstuff{X}$ acting on $\{\munit,\obstuff{X}\}$ is $\begin{psmallmatrix}0 & 1\\ 1& 1\end{psmallmatrix}$ has eigenvalues $\phi$ and $-(\phi-1)$. The largest eigenvalue $\phi$ gives the \textit{dominating} growth, and the absolute value of the second largest eigenvalue $\phi-1$ gives a bound for the variance.
\end{Example}

For $\catstuff{C}$ with finitely many isomorphism classes of indecomposable objects\emph{, called a finite growth problem,} 
there is a general answer to all questions in 
\autoref{Q:IntroductionMainQuestion}. The main player is the \emph{Perron--Frobenius (PF for short) dimension} which is the largest eigenvalue of the action matrix associated with the given growth problem.

For arbitrary additive Krull--Schmidt monoidal categories there is not much we can say about growth problems. However, we will add some natural assumptions (specified below) and address the following:
\begin{enumerate}[label=\emph{\upshape(\roman*)}]

\item Below we will attach $\pfdimplain\obstuff{X}\in\Rplus\cup\{\infty\}$ to $\obstuff{X}$ such that $(\pfdimplain\obstuff{X})^{n}$ is the dominating growth of $\bsymbol_{n}$. We show that this works in two cases 
using a generalization of the PF dimension which is the usual PF dimension for finite growth problems (and we write $\pfdimplain$ if the difference does not play a role):
\begin{enumerate}[label=\emph{\upshape(\arabic*)}]

\item Whenever a certain condition is satisfied which
we call \emph{PF admissible}. A key in this definition is that the associated PF dimension $\pfdim$ can be approximated using data from finite graphs (which is why we use $f$ as a subscript in the notation). 

\item Whenever a certain random walk associated to $\obstuff{X}$ is \emph{sustainably positively recurrent}. 
This case is a special case of the previous one, but more well-behaved so that we can say more; see (ii) up next.

\end{enumerate}
In particular, in these cases we have the \emph{exponential growth theorem}
$\lim_{n\to\infty}\sqrt[n]{\bsymbol_{n}}=\pfdim\obstuff{X}$.

\item If the category generated by $\obstuff{X}$ is sustainably positively recurrent,
then $\pfdimthree\obstuff{X}\neq\infty$ and 
we will give an explicit expression for $\asymbol(n)$ 
with $\bsymbol(n)\sim\asymbol(n)$. In this case we also 
determine the variance.

\end{enumerate}
In the main body of the paper we will define
our generalization of the PF dimension, \emph{(positive and null) recurrent} and \emph{transient categories} (and growth problems), following the classical notions in the theory of random walks, as well as PF admissible and 
sustainably positively recurrent growth problems. The rough mnemonic is as follows:
\begin{enumerate}[label=$\blacktriangleright$]

\item The (positive and null) recurrent growth problems, although not necessarily finite, behave like the finite case; in particular, if they are sustainably positively recurrent. This is the nicest possible case beyond finite growth problems.

\item The PF admissible growth problems generally behave quite differently from the finite case, but there is still a reasonable theory of PF dimensions.

\item For all other growth problems things become messy, and we do not know any general method for approaching such problems.

\end{enumerate}

\subsection{Key examples}

All cases where $\catstuff{C}$ has only finitely many indecomposable objects up to isomorphism are sustainably positively recurrent and PF admissible.
We discuss two additional classes of examples with potentially infinitely many isomorphism classes of indecomposable objects:
\begin{enumerate}

\item $(G,\obstuff{X})=\big(\rep(G),\obstuff{X}\big)$: Here $\rep(G)$ means finite dimensional $G$-representations over an arbitrary field, for $G$ a finite group and $\obstuff{X}\in\rep(G)$ arbitrary. This case is always sustainably positively recurrent. More generally, all growth problems in finite tensor categories are sustainably positively recurrent, and we will explore these as well. \textit{(Sustainably positively recurrent.)}

\item Again $(G,\obstuff{X})$: This time we consider $\rep(G)$ for $G$ a simple reductive group over a field of characteristic zero with a representation $\obstuff{X}\in\rep(G)$ in defining characteristic ({\ie} over the same field).
We show that this case is (almost always) transient, but PF admissible with $\pfdim\obstuff{X}\neq\infty$. \textit{(PF admissible.)}

\end{enumerate}
In the process, we also address (though some details are left to the reader):
\begin{enumerate}[resume]

\item $(\sbim,\obstuff{X})$: We discuss the category of Soergel bimodules $\sbim$ for an arbitrary finite Coxeter group $W$, 
and $\obstuff{X}$ a generating object, or any $\obstuff{X}$ such that some monoidal power of it has the indecomposable Soergel bimodule for the longest word of $W$ appearing as a direct summand (for $\sbim$ we use generating object to refer to either case). This case is sustainably positive recurrent. We will see that 
$\pfdim\obstuff{X}\neq\infty$ by giving an explicit formula for $\pfdim\obstuff{X}$ and for $\asymbol(n)$ with $\bsymbol(n)\sim\asymbol(n)$.
\textit{(Both, sustainably positively recurrent and PF admissible.)}

\item We also discuss examples where $\pfdim\obstuff{X}=\infty$.
Our examples include $(\mathrm{SL}_{\N}(\C),\C^{\N})$ and (generalized) Deligne categories. This gives new proofs that these categories have superexponential growth.

\end{enumerate}

Finally, in the process, we classify
for which $(G,\obstuff{X})$ is (positive and null) recurrent.

\subsection{Wrap up}

Before we get started, here are a few remarks.

\begin{Remark}\label{R:IntroHistory}
Our paper generalizes results that are known in the literature:
\begin{enumerate}

\item Our theorems applied to categories with finitely many indecomposable objects up to isomorphism recover the setting in \cite{LaTuVa-growth-pfdim}.

\item The exponential growth theorem 
partially generalizes \cite[Theorems 1.4 and 1.9]{CoOsTu-growth}.

\item Our classification of the (positive and null) recurrent
groups and representations is a version of Polya's random walk theorem \cite{Po-random-walks}.

\end{enumerate}
This will be explained in the main body of the paper.
\end{Remark}

\begin{Remark}\label{R:IntroRPlus}
All the above immediately generalize to certain algebras which have a basis with 
structure constants in $\Rplus$.
We will use this setting below. The connection to the above is 
that Grothendieck rings of additive Krull--Schmidt monoidal categories are $\Rplus$-algebras.
\end{Remark}

\begin{Remark}
Magma and Mathematica code for some 
of the calculations in this paper are available on GitHub: \cite{LaTuVa-pfdim-inf-code}.
\end{Remark}

\noindent\textbf{Acknowledgments.}
We are grateful to Kevin Coulembier, Pavel Etingof, David He, Eoghan McDowell, Jensen O'Sullivan, Victor Ostrik and A Referee for many helpful discussions about growth problems in general. We would also like to thank Nate Harman and Victor Ostrik for explaining the fusion rules of the Deligne and Delannoy categories to us. DT expresses particular appreciation for the enduring support encapsulated in the mantra, ``Definitions are more important than proofs.''

We utilized ChatGPT for proofreading support in the preparation of this paper.

DT was sponsored by the ARC Future Fellowship FT230100489, and PV was supported by the Fonds de la Recherche Scientifique-FNRS under Grant No. CDR-J.0189.23.


\section{The basic setting}\label{S:Basics}


The following mildly generalizes the usual notion of based $\Rplus$-algebras 
as {\eg} in \cite[Section 2]{KiMa-based-algebras}.
Let $\K$ be a unital subring of $\C$, e.g. $\K=\Z$. A \emph{based $\Rplus$-algebra} is a pair $(R,C)$ where $C=(c_{i})_{i\in I}$ with $1\in C$ is a $\K$-basis of a $\K$-algebra $R$ 
such that all structure constants lie in $\Rplus$, {\ie}:
\begin{gather}\label{Eq:BasicsSum}
c_{i}c_{j}=\sum_{k\in I}m_{ij}^{k}c_{k}\quad\text{with}\quad m_{ij}^{k}\in\Rplus.
\end{gather}
Note that $m_{ij}^{k}\in\K\cap\Rplus$, so that $m_{ij}^{k}\in\N$ when $\K=\Z$. We usually only write $R$ for $(R,C)$ if no confusion can arise. Note that the sum in \autoref{Eq:BasicsSum} is finite since $C$ is a $\K$-basis.

\begin{Remark}\label{R:GrowthRateCountable}
The reader who wants to consider $\Rplus$-algebras without the condition $1\in C$
can split their growth problem along an idempotent 
decomposition of the identity, so $1\in C$ 
is not a restriction.

By taking a subring of $R$ if necessary, we can, and will, always assume that $I$ is countable and that the identity is the element $c_{0}$ of the basis $C$. This is no restriction for growth problems as we will see below.
\end{Remark}

We denote by $\Rplus C$ the subset of $R$ of finite $\Rplus$-linear combinations of elements of the basis $C$ with nonnegative real coefficients, and similarly for $(\Rplus\cap\K)C$. Given $c\in(\Rplus\cap\K)C$, 
we write $cc_{j}=\sum_{k\in I}m_{c,j}^{k}c_{k}$ and
we denote $c^{n}=\sum_{i\in I}m_{n}^{i}(c)c_{i}$. We define the function
\begin{gather*}
\bsymbol^{R,c}\colon\N\to\Rplus,n\mapsto
\bsymbol(n)=\bsymbol^{R,c}(n):=\sum_{i\in I}m_{n}^{i}(c).
\end{gather*}
We are interested in the asymptotic behavior of this function. In particular, we use 
the terminology \emph{growth problems} associated to 
$(R,c)$ in the same way as in \autoref{Q:IntroductionMainQuestion}.

\begin{Convention}\label{N:GrowthRateProblems}
If not explicitly stated otherwise, representations are always finite dimensional.
\end{Convention}

\begin{Example}\label{E:GrowthRateProblems}
\emph{Good examples} of growth problems (coming from categories) that we can treat nicely are: 
\begin{enumerate}

\item Let $\catstuff{C}=\langle\obstuff{X}\rangle$ be the additive Krull--Schmidt monoidal category generated by $\obstuff{X}$ meaning that we take direct summands of finite sums of objects $\obstuff{X}^{\otimes n}$, where $n\in\N$ (we call such $\obstuff{X}$ generating objects). This example is good if $\catstuff{C}$ is a finite tensor category over an arbitrary field $\mathbb{F}$ in the sense of 
\cite[Definition 1.8.5 and Definition 4.1.1]{EtGeNiOs-tensor-categories} or finite with respect to isomorphism classes of indecomposable objects.
This includes:
\begin{enumerate}

\item $(H,\obstuff{X})=\big(\rep(H),\obstuff{X}\big)$
for $H$ a finite dimensional Hopf $\mathbb{F}$-algebra and $\obstuff{X}$ 
a finite dimensional $H$-representation over $\mathbb{F}$ that generates $\rep(H)$.
Specifically, $H$ could be the group ring of a finite group and $\obstuff{X}$ a faithful $H$-representation.

\item $(\sbim,\obstuff{X})$, for a generating object $\obstuff{X}$, as in \autoref{S:Intro}, see for example \cite{ElWi-soergel-calculus} for details about the category $\sbim$.

\end{enumerate}

\item $(G,\obstuff{X})=\big(\rep(G),\obstuff{X}\big)$ where $G$ is a simple reductive group and $\obstuff{X}$ is a $G$-representation, both in characteristic zero. (The assumption on $G$ being simple can be relaxed, see \autoref{R:GrowthRateMainTheoremGL}.)

\item $(\catstuff{C},\obstuff{X})$ for either of the following: $\catstuff{C}=\mathrm{SL}_{\N}(\C)$ and 
$\obstuff{X}=\C^{\N}$ the defining representation;
$\catstuff{C}=\rep(S_{t},\C)$ for $t\in\C\setminus\Z$ and 
$\obstuff{X}$ the generating object, following 
the conventions in \cite{CoOs-blocks-deligne-sym-group}, and $\catstuff{C}$ a Delannoy category with $\obstuff{X}$ the (defining) generating object following the conventions in \cite{HarSnoSny-delannoycat} (and \cite[Theorem 13.2]{HarSno-oligomorphic} shows that this is a good example).

\end{enumerate}
We will revisit these several times below.
\end{Example}

\begin{Remark}\label{R:GrowthRateGorthendieck}
In an additive Krull--Schmidt monoidal category $\catstuff{C}$ the \emph{fusion graph}
of an object $\obstuff{X}$ is the graph with vertices 
corresponding to the indecomposable objects in $\catstuff{C}$
and outgoing edges of some vertex $\obstuff{Y}$ 
going to the indecomposable direct summands of $\obstuff{X}\otimes\obstuff{Y}$.
The below is motivated by this terminology.
\end{Remark}

Let us consider the \emph{(oriented and weighted) fusion graph} $\Gamma=\Gamma(c)$ 
associated to some growth problem $(R,c)$. This graph is defined as follows:
\begin{enumerate}[label=\emph{\upshape(\roman*)}]

\item The vertices of $\Gamma$ correspond to the basis elements $c_{i}\in C$ appearing 
in $c^{n}$ for some $n$, {\ie} those $c_{i}$ for which some $m_{n}^{i}\neq 0$.

\item There is an edge of labeled $m_{c,j}^{i}$ from the vertex $c_{j}$ to the vertex $c_{i}$.

\end{enumerate}
We identify the basis elements $c_{i}\in C$ with the vertices of $\Gamma$.
After fixing some ordering of the vertices, the matrix associated to a fusion graph 
$\Gamma$ is called the \emph{action matrix} $M(\Gamma)=(m_{c,j}^{k})_{k,j}$.

\begin{Example}\label{E:GrowthRateSLTwoOne}
The key example the reader should keep in mind is the growth 
problem for $\big(\mathrm{SL}_{2}(\C),\C^{2}\big)$. In this case the fusion graph is $\N$:
\begin{gather}\label{Eq:GrowthRateSLTwoOne}
\begin{gathered}
\Gamma(\mathrm{SL}_{2})=
\begin{tikzpicture}[anchorbase]
\begin{scope}[every node/.style={circle,inner sep=0pt,text width=6mm,align=center,draw=white,fill=white}]
\node (A) at (0,0) {$\bullet$};
\node (B) at (2,0) {$\bullet$};
\node (C) at (4,0) {$\bullet$};
\node (D) at (6,0) {$\bullet$};
\node (E) at (8,0) {$\bullet$};
\node (F) at (10,0) {$\bullet$};
\end{scope}
\node (G) at (12,0) {$\dots$};
\begin{scope}[>={Stealth[black]},
every edge/.style={draw=orchid,very thick}]
\path [->] ($(A)+(0.25,0.075)$) edge ($(B)+(-0.25,0.075)$);
\path [->] ($(B)+(0.25,0.075)$) edge ($(C)+(-0.25,0.075)$);
\path [->] ($(C)+(0.25,0.075)$) edge ($(D)+(-0.25,0.075)$);
\path [->] ($(D)+(0.25,0.075)$) edge ($(E)+(-0.25,0.075)$);
\path [->] ($(E)+(0.25,0.075)$) edge ($(F)+(-0.25,0.075)$);
\path [->] ($(F)+(0.25,0.075)$) edge ($(G)+(-0.25,0.075)$);
\path [<-] ($(A)+(0.25,-0.075)$) edge ($(B)+(-0.25,-0.075)$);
\path [<-] ($(B)+(0.25,-0.075)$) edge ($(C)+(-0.25,-0.075)$);
\path [<-] ($(C)+(0.25,-0.075)$) edge ($(D)+(-0.25,-0.075)$);
\path [<-] ($(D)+(0.25,-0.075)$) edge ($(E)+(-0.25,-0.075)$);
\path [<-] ($(E)+(0.25,-0.075)$) edge ($(F)+(-0.25,-0.075)$);
\path [<-] ($(F)+(0.25,-0.075)$) edge ($(G)+(-0.25,-0.075)$);
\end{scope}
\end{tikzpicture}
,
\\
M(\mathrm{SL}_{2})=
\begin{psmallmatrix}
0 & 1 & 0 & 0 & 0 & \dots
\\
1 & 0 & 1 & 0 & 0 & \dots
\\
0 & 1 & 0 & 1 & 0 & \dots
\\
0 & 0 & 1 & 0 & 1 & \dots
\\
0 & 0 & 0 & 1 & 0 & \dots
\\
\vdots & \vdots & \vdots & \vdots & \vdots & \vdots
\end{psmallmatrix}.
\end{gathered}
\end{gather}
We will come back to this growth problem several times throughout the paper.
\end{Example}

\begin{Example}\label{E:GrowthRateGrothendieck}
The Grothendieck group of an additive Krull--Schmidt monoidal category $\catstuff{C}$ is an example of a based $\Rplus$-algebra with $\K=\Z$ with basis given by 
the classes of indecomposable objects. 
The growth problems as in 
\autoref{Q:IntroductionMainQuestion} can then be reformulated in terms of 
growth problems for based $\Rplus$-algebras, and we do this silently throughout.
In this case the fusion graph $\Gamma(\obstuff{X})$ for $\obstuff{X}\in\catstuff{C}$
is the fusion graph on the subcategory generated by $\obstuff{X}^{\otimes n}$ in the usual sense. 
\end{Example}

\begin{Lemma}\label{L:GrowthRateBasicFusion}
For a fusion graph $\Gamma$ for a growth problem we have:
\begin{enumerate}

\item $\Gamma$ has countably many vertices.

\item There is a path from $c_{0}$ to any other 
vertex of $\Gamma$. In particular, 
$\Gamma$ is connected in the unoriented sense.

\item Every vertex of $\Gamma$ is of 
finite degree.

\end{enumerate}
\end{Lemma}

\begin{proof}
This follows because we only consider $c^{n}$, which, for fixed $n$, has only finitely many nonzero coefficients in terms of $C$.
\end{proof}


\section{The dominating growth via truncations}\label{S:GrowthRate}


The \emph{naive cutoff} of $\Gamma$ is
the sequence $(\Gamma_{k})_{k\in\N}$ of subgraphs 
of $\Gamma$ of the form
\begin{gather*}
(\Gamma_{k})_{k\in\N}=\big(\Gamma_{0}=\{c_0\}\subset\Gamma_{1}\subset\Gamma_{2}\subset\dots\subset\Gamma\big)
\text{ with }\bigcup_{k\in\N}\Gamma_{k}=\Gamma,
\end{gather*}
such that $\Gamma_{k}$ is the induced subgraph of $\Gamma$ 
that has precisely the vertices that 
can be reached with a path of length $\leq k$ from $c_{0}$.
Note that \autoref{L:GrowthRateBasicFusion}.(c) implies that $\Gamma_{k}$ will be finite.

\begin{Definition}\label{D:GoodFiltration}
A \emph{good filtration} of a growth problem 
$(R,C)$ 
is a sequence of (not necessarily associative) $\Rplus$-algebras $\big((R_{k},C_{k})\big)_{k\in\N}$ over $\K$ such that:
\begin{enumerate}

\item $C_{k}\subset C$ is finite and
\begin{gather*}
(C_{k})_{k\in\N}=
\big(C_{0}=\{c_0\}\subset C_{1}\subset\dots\subset C_{k}\subset\dots\subset C\big)\text{ with }\bigcup_{k\in\N}C_{k}=C.
\end{gather*}

\item $c_{0}$ is the unit of all $R_{k}$.

\item The multiplication in $R_{k}$ is given by taking the product in $(R,C)$ followed by the projection to $(R_{k},C_{k})$.

\end{enumerate}
The \emph{naive filtration} of $(R,C)$ is the good filtration where $C_{k}$ are the vertices of $\Gamma_{k}$ for the naive cutoff.
\end{Definition}

\begin{Example}\label{E:GrowthRateSLTwoOneCutOff}
In \autoref{Eq:GrowthRateSLTwoOne} the naive cutoff consists of 
$\Gamma_{k}$ being a line graph with $k$ vertices. In 
\autoref{Eq:GrowthRateMainTheoremKlein} below, all the $\Gamma_{k}$ for $k>1$ will include the vertex labeled $\prin[4]$. 
\end{Example}

For a good filtration, in $(R_{k},C_{k})$ (this is of finite rank over $\K$) we can define $\pfdimk c_{i}$ to be the \emph{PF eigenvalue} of the action matrix of $c_{i}$ with respect to $C_{k}$.

\begin{Lemma}\label{L:GrowthLeq}
We have $\pfdimk c_{i}\leq\pfdimk[k+1]c_{i}$ and
$\lim_{k\to\infty}\pfdimk c_{i}\in\Rplus\cup\{\infty\}$ is well-defined.
\end{Lemma}

\begin{proof}
The first statement follows from the fact that subgraphs of a graphs cannot have larger PF eigenvalues than the original graph. The second statement then follows since the previous point gives us an increasing sequence.
\end{proof}

\begin{Definition}\label{D:GrowthRatePFDim}
Given $c\in(\Rplus\cap\K)C$, define the (good filtration) \emph{PF dimension} as the limit $\pfdim c=\lim_{k\to\infty}\pfdimk c\in\Rplus\cup\{\infty\}$. 
\end{Definition}

Note that \autoref{D:GrowthRatePFDim} depends on the choice of a good filtration, see \autoref{R:GrowthRateMainTheoremGL} below for an almost example (illustrating why one needs to be careful with the choice of the filtration).

\begin{Lemma}\label{L:GrowthRatePFFinite}
If $|I|<\infty$, then 
\autoref{D:GrowthRatePFDim} agrees with the usual 
definition of PF dimension using the largest eigenvalue of the associated matrix (or, equivalently, graph).
\end{Lemma}

\begin{proof}
Immediate.
\end{proof}

\begin{Definition}\label{D:GrowthRatePFAdmissible}
A growth problem $(R,C)$ is called \emph{PF admissible} if there exists a good filtration such that:
\begin{enumerate}

\item $\pfdim$ is superadditive and subpowermultiplicative:
\begin{gather*}
\pfdim c+\pfdim d\leq\pfdim(c+d),
\\
\pfdim c^{n}\leq(\pfdim c)^{n},
\end{gather*}
for all $c,d\in(\Rplus\cap\K)C$, $n\in\N$.

\item $\pfdim c_{i}\geq 1$ for all $c_{i}\in C$.

\end{enumerate}
We always associate any such good filtration to a PF admissible growth problem.
We use the same terminology for the fusion graphs and action matrices themselves.
\end{Definition}

\begin{Remark}
If $R$ is of finite $\K$-rank, then two good filtrations will agree on all but finitely many pages of the filtrations.
Abusing terminology a bit, we do not need to and will not specify a specific good filtration.
\end{Remark}

\begin{Example}\label{E:GrowthRatePFAdmissibleFusion}
A growth problem coming from a transitive unital $\Zplus$-ring 
of finite $\Z$-rank is PF admissible. (This includes decategorifications of fusion categories.)
The condition $\pfdim c_{i}\geq 1$ for all $c_{i}\in C$ holds in this case, see e.g. \cite[Proposition 3.3.4]{EtGeNiOs-tensor-categories}. 
In fact, this assumption is motivated from these rings.
\end{Example}

\begin{Example}\label{E:GrowthRatePFAdmissible}
By a classical result of Kronecker, as summarized in e.g. \cite[Lemma 3.3.14]{EtGeNiOs-tensor-categories}, every growth problem corresponding to a finite graph satisfies \autoref{D:GrowthRatePFAdmissible}.(b) unless the PF dimension is zero.
\end{Example}

\begin{Notation}\label{N:GrowthRateMainTheorem}
We use the usual \emph{Bachmann--Landau} (also called Landau--Bachmann or just Landau or even capital O) \emph{notation}:
A function $f\colon\N\to\Rplus$ satisfies $f\in\Theta(g)$ if there exist constants $A,B\in\R_{>0}$ such that $A\cdot g(n)\leq f(n)\leq B\cdot g(n)$ for 
all $n_{0}<n$ for some fixed $n_{0}\in\N$. If one has the upper bound, then 
we write $f\in O(g)$, and if the lower bound holds, then we write $f\in\Omega(g)$. 

Moreover, a function $f$ is of \emph{superexponential growth} if
$f\in\Omega(\gamma^{n})$ for all $\gamma\in\R_{>1}$ (in particular, $\sqrt[n]{f(n)}$ is unbounded in this case),
and is of \emph{subexponential growth} if $f\in O(\gamma^{n})$
for all $\gamma\in\R_{>1}$.
\end{Notation}

\begin{Theorem}\label{T:GrowthRateMainTheorem}
Let $(R,c)$ be a growth problem.
\begin{enumerate}

\item For a good filtration we have
\begin{gather}\label{Eq:GrowthRateMainTheoremOne}
\bsymbol(n)\in\Omega(\gamma^{n})
\text{ for all }\gamma<\pfdim c.
\end{gather}

\end{enumerate}

Assume now that $(R,c)$ is PF admissible.

\begin{enumerate}[resume]

\item If $\pfdim c\neq\infty$, then
\begin{gather}\label{Eq:GrowthRateMainTheoremTwo}
\bsymbol(n)\in O\big((\pfdim c)^{n}\big).
\end{gather}

\item We have the \emph{exponential growth theorem}:
\begin{gather*}
\beta=\lim_{n\to\infty}\sqrt[n]{\bsymbol_{n}}=\pfdim c.
\end{gather*}

\end{enumerate}
In particular, $\bsymbol_{n}$ grows superexponentially if and only if
$\pfdim c=\infty$.
\end{Theorem}

In general the PF dimension depends on the choice of a good filtration. However, if one has two different good filtrations such that the growth problem is PF admissible, then the PF dimension is the same for both of them by \autoref{T:GrowthRateMainTheorem}.(c).

\begin{proof}[Proof of \autoref{T:GrowthRateMainTheorem}]
As in the proof of \cite[Theorem 3B.2]{LaTuVa-growth-pfdim},
the number $\bsymbol_{n}$ can be computed as the column sum 
of the matrix $M(\Gamma)^{n}$ for the first column corresponding 
to the identity. (We will use this throughout.) To see this we consider 
the equation
\begin{gather*}
M(\Gamma)c(n-1)=c(n)
\xrightarrow{\text{iterate}}
M(\Gamma)^{n}c(0)=c(n),
\end{gather*}
where $c(k)=\big(c_{0}(k),c_{1}(k),\dots\big)\in\Rplus^{\N}$ are 
vectors such that their $i$th entry is the 
multiplicity of $c_{i}$ in $c^{k}$, and $c(0)=(1,0,\dots)^{T}$ with the 
one is in the slot of $c_{0}=1$. Observe $M(\Gamma)^{n}c(0)$ is the first column of $M(\Gamma)^{n}$.

\textit{(a).} By an appropriate version of Feteke's subadditive lemma, the sequence $(\sqrt[n]{\bsymbol_{n}})_{n\in\N}$ has a limit
\begin{gather*}
\lim_{n\to\infty}\sqrt[n]{\bsymbol_{n}}=\beta\in\Rplus\cup\{\infty\}.
\end{gather*}
Similarly, if we denote by $(\bsymbol^{(k)}_{n})_{n\in\N}$ the respective growth problem in the cutoff $R_{k}$ for $c\in R_{k}$ and set $\bsymbol^{(k)}_{n}=0$ whenever $c\notin R_{k}$ (this happens for only finitely many $k\in\N$), then
\begin{gather*}
\lim_{n\to\infty}\sqrt[n]{\bsymbol^{(k)}_{n}}=\beta_{k}\in\Rplus\cup\{\infty\},
\end{gather*} 
exists. Note hereby that $\beta_{k}=\pfdimk c$, by 
\cite[Theorem 1]{LaTuVa-growth-pfdim} (which uses PF theory of finite graphs).
Moreover, the definition of a good filtration \autoref{D:GoodFiltration} implies that $\bsymbol^{(k)}_{n}\leq\bsymbol_{n}$ which in turn implies that $\beta_{k}\leq\beta$ for all $k,n\in\N$. Therefore, $\pfdim c\leq \beta$ by \autoref{D:GrowthRatePFDim}. Hence, 
$\bsymbol(n)\in\Omega(\gamma^{n})$
for all $\gamma<\pfdim c$.

\textit{(b).}  We go back to the definition of $\bsymbol_{n}$. We have $c^{n}=\sum_{i\in I}m_{n}^{i}(c)c_{i}$ and we obtain:
\begin{gather*}
(\pfdim c)^{n}\geq\pfdim c^{n}\geq\sum_{i\in I}m_{n}^{i}(c)\pfdim c_{i}\geq\sum_{i\in I}m_{n}^{i}(c)=\bsymbol_{n},
\end{gather*}
the inequalities following from the PF admissibility condition.

\textit{(c).} This follows from the just established 
\autoref{Eq:GrowthRateMainTheoremOne} and \autoref{Eq:GrowthRateMainTheoremTwo}.
\end{proof}

\begin{Example}\label{E:GrowthRateSLTwoTwo}
Coming back to \autoref{E:GrowthRateSLTwoOne}, the cutoff $(\Gamma_{i})_{i\in\N}$ of $\Gamma(\mathrm{SL}_{2})$ is such that $\Gamma_{i}$ is the finite subgraph on the first $i$ vertices exemplified by:
\begin{gather*}
\Gamma_{7}
=
\begin{tikzpicture}[anchorbase]
\begin{scope}[every node/.style={circle,inner sep=0pt,text width=6mm,align=center,draw=white,fill=white}]
\node (A) at (0,0) {$\bullet$};
\node (B) at (2,0) {$\bullet$};
\node (C) at (4,0) {$\bullet$};
\node (D) at (6,0) {$\bullet$};
\node (E) at (8,0) {$\bullet$};
\node (F) at (10,0) {$\bullet$};
\node (G) at (12,0) {$\bullet$};
\end{scope}
\begin{scope}[>={Stealth[black]},
every edge/.style={draw=orchid,very thick}]
\path [->] ($(A)+(0.25,0.075)$) edge ($(B)+(-0.25,0.075)$);
\path [->] ($(B)+(0.25,0.075)$) edge ($(C)+(-0.25,0.075)$);
\path [->] ($(C)+(0.25,0.075)$) edge ($(D)+(-0.25,0.075)$);
\path [->] ($(D)+(0.25,0.075)$) edge ($(E)+(-0.25,0.075)$);
\path [->] ($(E)+(0.25,0.075)$) edge ($(F)+(-0.25,0.075)$);
\path [->] ($(F)+(0.25,0.075)$) edge ($(G)+(-0.25,0.075)$);
\path [<-] ($(A)+(0.25,-0.075)$) edge ($(B)+(-0.25,-0.075)$);
\path [<-] ($(B)+(0.25,-0.075)$) edge ($(C)+(-0.25,-0.075)$);
\path [<-] ($(C)+(0.25,-0.075)$) edge ($(D)+(-0.25,-0.075)$);
\path [<-] ($(D)+(0.25,-0.075)$) edge ($(E)+(-0.25,-0.075)$);
\path [<-] ($(E)+(0.25,-0.075)$) edge ($(F)+(-0.25,-0.075)$);
\path [<-] ($(F)+(0.25,-0.075)$) edge ($(G)+(-0.25,-0.075)$);
\end{scope}
\end{tikzpicture}
.
\end{gather*}
We will prove later that this growth problem with the naive truncation is PF admissible, see \autoref{P:GrowthRatePFAdmissible}.
It is easy to see that $\pfdimk[i]\Gamma_{i}=2\cos\big(\pi/(i+1)\big)$.
Then 
$\pfdim\Gamma(\mathrm{SL}_{2})=
\lim_{i\to\infty}2\cos\big(\pi/(i+1)\big)=2$. 
Thus, \autoref{T:GrowthRateMainTheorem} implies that the dominating growth factor of $\bsymbol_{n}$ is $2^{n}$.
Moreover, dividing by the dominating growth factor, the log plots
\begin{gather*}
\begin{tikzpicture}[anchorbase]
\node at (0,0) {\includegraphics[height=4.6cm]{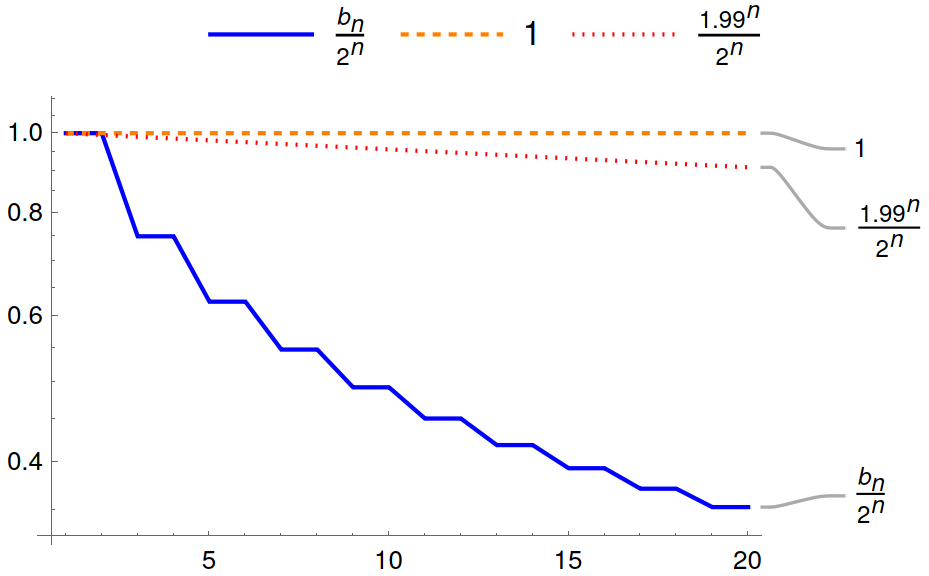}};
\end{tikzpicture}
,\quad
\begin{tikzpicture}[anchorbase]
\node at (0,0) {\includegraphics[height=4.6cm]{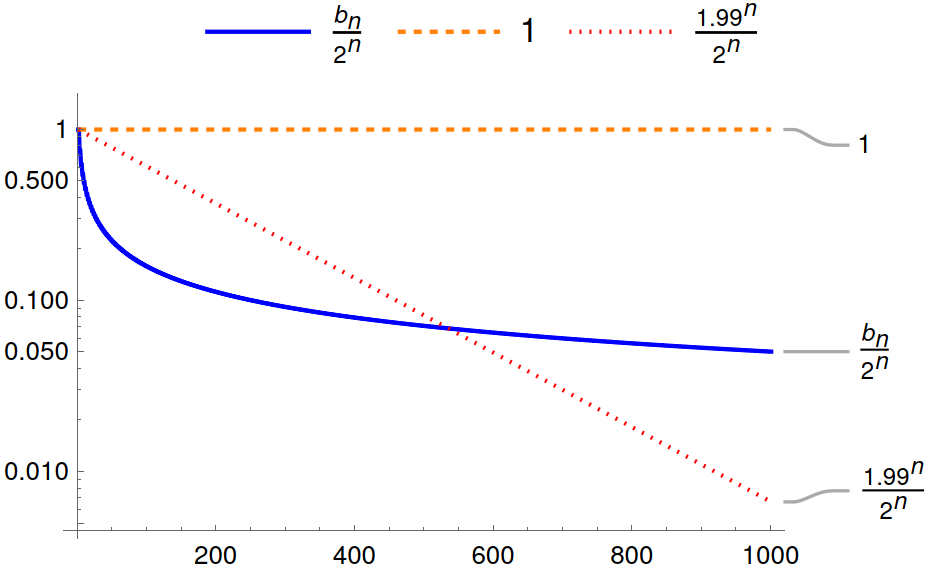}};
\end{tikzpicture}
,
\end{gather*}
shows \autoref{Eq:GrowthRateMainTheoremOne} for $\gamma=1.99$ and \autoref{Eq:GrowthRateMainTheoremTwo}.
\end{Example}

Next, a ``nonexample'':

\begin{Example}\label{E:GrowthRateMainTheoremNecessary}
We give an example where the PF admissibility condition is necessary for \autoref{T:GrowthRateMainTheorem}. Suppose that our $\Rplus$-algebra has a basis $(c_{i})_{i\in\N}$ with multiplication $c_{i}c_{j}=c_{i+j}$ and that $c=\alpha c_{0}+c_{1}$ for some $\alpha\in\Rplus$. Then $c$ satisfies $cc_{i}=\alpha\cdot c_{i}+c_{i+1}$ for all $i\in\N$. The action matrix is an infinite Jordan block with $\alpha$-labeled loops and we have
\begin{gather*}
\Gamma(Jo)=
\raisebox{-0.22cm}{$\begin{tikzpicture}
\begin{scope}[every node/.style={circle,inner sep=0pt,text width=6mm,align=center,draw=white,fill=white}]
\node (A) at (0,0) {$\bullet$};
\node (B) at (2,0) {$\bullet$};
\node (C) at (4,0) {$\bullet$};
\node (D) at (6,0) {$\bullet$};
\node (E) at (8,0) {$\bullet$};
\node (F) at (10,0) {$\bullet$};
\end{scope}
\node (G) at (12,0) {$\dots$};
\begin{scope}[>={Stealth[black]},
every edge/.style={draw=orchid,very thick}]
\path [->] (A) edge (B);
\path [->] (B) edge (C);
\path [->] (C) edge (D);
\path [->] (D) edge (E);
\path [->] (E) edge (F);
\path [->] (F) edge (G);
\path [->] (A) edge [loop above] node{$\alpha$} (A);
\path [->] (B) edge [loop above] node{$\alpha$} (B);
\path [->] (C) edge [loop above] node{$\alpha$} (C);
\path [->] (D) edge [loop above] node{$\alpha$} (D);
\path [->] (E) edge [loop above] node{$\alpha$} (E);
\path [->] (F) edge [loop above] node{$\alpha$} (F);
\end{scope}
\end{tikzpicture}$}
,
\\
M(Jo)=
\begin{psmallmatrix}
\alpha & 0 & 0 & 0 & 0 & \dots
\\
1 & \alpha & 0 & 0 & 0 & \dots
\\
0 & 1 & \alpha & 0 & 0 & \dots
\\
0 & 0 & 1 & \alpha & 0 & \dots
\\
0 & 0 & 0 & 1 & \alpha & \dots
\\
\vdots & \vdots & \vdots & \vdots & \vdots & \vdots
\end{psmallmatrix}
,
\end{gather*}
as the fusion graph. The naive cutoff is the sequence of finite Jordan blocks, in particular, $\pfdim c_{i}=0$ for $i>0$ and $\alpha=0$, so this example does not satisfy \autoref{D:GrowthRatePFAdmissible}.(b).
And in fact we get:
\begin{gather*}
c^{n}=\sum_{i=0}^{n}\binom{n}{i}\alpha^{n-i}\cdot c_{i},
\end{gather*}
and therefore $\bsymbol_{n}=(1+\alpha)^{n}$ but $\pfdim c=\alpha$.
\end{Example}

Here are key examples of categories with infinitely many indecomposable objects up to isomorphism that are PF admissible:

\begin{Proposition}\label{P:GrowthRatePFAdmissible}
Any growth problem in \autoref{E:GrowthRateProblems}.(b) is 
PF admissible.
\end{Proposition}

\begin{proof}
Recall that we only consider simple reductive groups.
We will make use of the \emph{Verlinde categories}
$\catstuff{V}_{k}(G)$ of level $k$, see 
for example \cite{Sa-quantum-roots-of-unity} for a nice summary of the main properties we need (when we refer to it we use the arXiv numbering). For example, that the PF dimension is the categorical dimension.

In these semisimple categories the indecomposable (= simple) objects
$\overline{L}(\lambda,k)$ 
are indexed by $\Lambda^{k}$, which are cutoffs of the dominant integral Weyl chamber, 
see \cite[Section 4]{Sa-quantum-roots-of-unity} for details. Moreover, we have
\begin{gather*}
\pfdimk\overline{L}(\lambda,k)=\dim_{\catstuff{V}_{k}(G)}\overline{L}(\lambda,k)=
\prod_{\beta>0}\frac{q^{\langle\lambda+\rho, \beta\rangle}-q^{-\langle\lambda+\rho,\beta\rangle}}{q^{\langle\rho,\beta\rangle}-q^{-\langle\rho,\beta\rangle}}
\end{gather*}
(the notation is as in \cite[Section 1]{Sa-quantum-roots-of-unity})
the so-called \emph{quantum Weyl dimension formula}. Here $q=q(k)$ is a certain complex root of unity.

The only thing we need to know about this formula is that
\begin{gather*}
\lim_{k\to\infty}\pfdimk\overline{L}(\lambda,k)
=
\prod_{\beta>0}\frac{\langle\lambda+\rho, \beta\rangle}{\langle\rho,\beta\rangle},
\end{gather*}
which is now the \emph{classical} Weyl dimension formula. This is well-known, since 
$k\to\infty$ corresponds to $q\to 1$.

Let us choose the cutoff such that the vertices
of $\Gamma_{k}$ can be 
matched with $\Lambda^{k}$. Using this choice 
the induced product matches the product in the Verlinde category, by the fusion rules as recalled in \cite[Section 5]{Sa-quantum-roots-of-unity}. (It is remarkable that the Verlinde product is a truncation.)
Hence, this shows that \autoref{D:GrowthRatePFAdmissible}.(b) holds.

Moreover, since the Verlinde categories are braided,
it follows that the action matrices $M(c,k)$ of any $\Rplus$-linear combination of the $[\overline{L}(\lambda,k)]$ commute with one another, and these matrices can be simultaneously diagonalized. It is easily checked that the vector
\begin{gather*}
\sum_{\lambda\in\Lambda^{k}}\pfdimk\overline{L}(\lambda,k)\cdot[\overline{L}(\lambda,k)]
\end{gather*}
is an eigenvector of $M(c,k)$. The classical PF theorem implies that the associated eigenvalue of $M(c,k)$ is $\pfdimk M(c,k)$.
It follows that we have subpowermultiplicativity, or more precisely even:
\begin{gather*}
\pfdimk\big(M(c,k)M(c^{\prime},k)\big)
=\pfdimk\big(M(c,k)\big)\pfdimk\big(M(c^{\prime},k)\big).
\end{gather*}
and superadditivity, ore more precisely even:
\begin{gather*}
\pfdimk\big(M(c,k)+M(c^{\prime},k)\big)
=\pfdimk\big(M(c,k)\big)+\pfdimk\big(M(c^{\prime},k)\big).
\end{gather*}
This then implies \autoref{D:GrowthRatePFAdmissible}.(a).
\end{proof}

\begin{Example}\label{R:GrowthRateMainTheoremGL}
The argument in the proof of \autoref{P:GrowthRatePFAdmissible}
also works, for example, for the nonsimple group
$\mathrm{GL}_{n}(\C)$ and its defining representation. However, the setting we use needs to be adjusted since the induced product is not the product in the Verlinde cutoff.

Explicitly, take $n=2$. Let us denote the vector representation of $\mathrm{GL}_{2}(\C)$ by $L(1,0)$. Then there exists simple $\mathrm{GL}_{2}(\C)$-representations $L(a,b)$ with $a,b\in\N,a\geq b$ such that $L(1,0)\otimes L(a,b)\cong L(a+1,b)\oplus L(a,b+1)$, with the convention that $L(a,b)$ is zero if $a,b\in\N,a\geq b$ is not satisfied. The fusion graph is then of the form (with the vertices on a grid using $(a,b)$ as coordinates):
\begin{gather*}
\begin{gathered}
	\Gamma(\mathrm{GL}_{2})=
	\begin{tikzpicture}[anchorbase]
		\begin{scope}[every node/.style={circle,inner sep=0pt,text width=6mm,align=center,draw=white,fill=white}]
			\node (A) at (0,0) {$\bullet$};
			\node (B) at (2,0) {$\bullet$};
			\node (C) at (4,0) {$\bullet$};
			\node (D) at (6,0) {$\bullet$};
			\node (E) at (8,0) {$\bullet$};
			\node (F) at (10,0) {$\bullet$};
			\node (B2) at (2,2) {$\bullet$};
			\node (C2) at (4,2) {$\bullet$};
			\node (D2) at (6,2) {$\bullet$};
			\node (E2) at (8,2) {$\bullet$};
			\node (F2) at (10,2) {$\bullet$};
			\node (C3) at (4,4) {$\bullet$};
			\node (D3) at (6,4) {$\bullet$};
			\node (E3) at (8,4) {$\bullet$};
			\node (F3) at (10,4) {$\bullet$};
			\node (D4) at (6,6) {$\dots$};
			\node (E4) at (8,6) {$\dots$};
			\node (F4) at (10,6) {$\dots$};
		\end{scope}
		\node (G) at (12,0) {$\dots$};
		\node (G2) at (12,2) {$\dots$};
		\node (G3) at (12,4) {$\dots$};
		\node (G4) at (12,6) {$\dots$};
		\begin{scope}[>={Stealth[black]},
			every edge/.style={draw=orchid,very thick}]
			\path [->] (A) edge (B);
			\path [->] (B) edge (C);
			\path [->] (C) edge (D);
			\path [->] (D) edge (E);
			\path [->] (E) edge (F);
			\path [->] (F) edge (G);
			\path [->] (B2) edge (C2);
			\path [->] (C2) edge (D2);
			\path [->] (D2) edge (E2);
			\path [->] (E2) edge (F2);
			\path [->] (F2) edge (G2);
			\path [->] (C3) edge (D3);
			\path [->] (D3) edge (E3);
			\path [->] (E3) edge (F3);
			\path [->] (F3) edge (G3);
			\path [->] (B) edge (B2);
			\path [->] (C) edge (C2);
			\path [->] (D) edge (D2);
			\path [->] (E) edge (E2);
			\path [->] (F) edge (F2);
			\path [->] (C2) edge (C3);
			\path [->] (D2) edge (D3);
			\path [->] (E2) edge (E3);
			\path [->] (F2) edge (F3);
			\path [->] (D3) edge (D4);
			\path [->] (E3) edge (E4);
			\path [->] (F3) edge (F4);
		\end{scope}
	\end{tikzpicture}
	,
	\\
	M(\mathrm{GL}_{2})=
    \begin{psmallmatrix}
0 & 0 & 0 & 0 & 0 & \dots
\\
J_{1} & 0 & 0 & 0 & 0 & \dots
\\
0 & J_{2}^{\prime} & 0 & 0 & 0 & \dots
\\
0 & 0 & J_{2} & 0 & 0 & \dots
\\
0 & 0 & 0 & J_{3}^{\prime} & 0 & \dots
\\
\vdots & \vdots & \vdots & \vdots & \vdots & \vdots
\end{psmallmatrix}
,
\\
\text{ with }
\begin{cases}
J_{k}=\text{ (lower triangular) Jordan block of size $k$ with eigenvalue $1$},
\\
J_{k}^{\prime}=J_{k}\text{ without the last column}.
\end{cases}
\end{gathered}
\end{gather*}
The naive filtration therefore gives $\pfdim\Gamma(\mathrm{GL}_{2})=0$ since all cutoffs are described by nilpotent matrices, and \autoref{D:GrowthRatePFAdmissible}.(b) is not satisfied. In contrast, the Verlinde filtration gives $\pfdim\Gamma(\mathrm{GL}_{2})=2$ (as the corresponding graph has backward arrows).
\end{Example}

\begin{Remark}\label{R:GrowthRateMainTheoremGen}
\autoref{P:GrowthRatePFAdmissible} and its proof also generalize \cite[Theorem 1.2]{LiShYaZh-pf-stuff}.
\end{Remark}



%

\begin{Definition}\label{D:GrowthRateSuper}
We say an $\Rplus$-algebra is of \emph{superexponential growth} if there exists some $c$ such that $b_{n}$ grows superexponentially.

We say an additive 
Krull--Schmidt monoidal category 
is of \emph{superexponential growth} 
if its Grothendieck ring is of superexponential growth 
in the above sense. 
\end{Definition}

\begin{Remark}\label{R:GrowthRateSuper}
\autoref{D:GrowthRateSuper} runs
in parallel to the definition of 
superexponential growth for abelian monoidal categories using 
the length of objects. However, note that:
\begin{enumerate}

\item An abelian growth problem could be of superexponential growth while its version defined using $b_{n}$ might not grow superexponentially. We, however, do not know any example where this happens.

\item If $b_{n}$ grows superexponentially, then
the abelian one does so as well.

\end{enumerate}
Hence, the notion from \autoref{D:GrowthRateSuper} ``grows (potentially) slower'' than the abelian one.
\end{Remark}

The next example is simple but crucial to determine categories of superexponential growth.

\begin{Example}\label{E:GrowthRateStarGraph}
Consider the following graph:
\begin{gather*}
\Gamma(\star_{N})=
\begin{tikzpicture}[anchorbase]
\begin{scope}[every node/.style={circle,inner sep=0pt,text width=4mm,align=center,draw=black,fill=white}]
\node (Z) at (0,0) {$0$};
\node (A) at (-1,0) {$1$};
\node (B) at (0,1) {$2$};
\node (D) at (0,-1) {$N$};
\end{scope}
\node (C) at (1,0) {$\dots$};
\node at (0.5,0.5) {$\ddots$};
\node at (0.5,-0.5) {$\rddots$};
\begin{scope}[>={Stealth[black]},
every edge/.style={draw=orchid,very thick}]
\path [->] ($(Z)+(-0.25,0.075)$) edge ($(A)+(0.25,0.075)$);
\path [->] ($(Z)+(0.075,0.25)$) edge ($(B)+(0.075,-0.25)$);
\path [->] ($(Z)+(0.25,0.075)$) edge ($(C)+(-0.25,0.075)$);
\path [->] ($(Z)+(0.075,-0.25)$) edge ($(D)+(0.075,0.25)$);
\path [<-] ($(Z)+(-0.25,-0.075)$) edge ($(A)+(0.25,-0.075)$);
\path [<-] ($(Z)+(-0.075,0.25)$) edge ($(B)+(-0.075,-0.25)$);
\path [<-] ($(Z)+(0.25,-0.075)$) edge ($(C)+(-0.25,-0.075)$);
\path [<-] ($(Z)+(-0.075,-0.25)$) edge ($(D)+(-0.075,0.25)$);
\end{scope}
\end{tikzpicture}
,\quad
M(\star_{N})=
\begin{psmallmatrix}
0 & 1 & 1 & \dots & 1 \\
1 & 0 & 0 & \dots & 0 \\
1 & 0 & 0 & \dots & 0 \\
\vdots & \vdots & \vdots & \ddots & \vdots \\
1 & 0 & 0 & \dots & 0
\end{psmallmatrix}
.
\end{gather*}
This is called the \emph{star graph} with $N+1$ vertices.
It is easy to see that $\pfdimplain\Gamma(\star_{N})=\sqrt{N}$.

Let $\Gamma(Y)$ be the Young lattice, considered as an oriented graph by 
putting an orientation in both directions. The first few layers of this graph are:
\begin{gather*}
\Gamma(Y)\leftrightsquigarrow	
\scalebox{0.6}{$\begin{tikzpicture}[>=latex,line join=bevel,scale=0.7,anchorbase]
\node (node_0) at (148.0bp,7.5bp) [draw,draw=none] {${\emptyset}$};
\node (node_1) at (148.0bp,60.5bp) [draw,draw=none] {${\def\lr#1{\multicolumn{1}{|@{\hspace{.6ex}}c@{\hspace{.6ex}}|}{\raisebox{-.3ex}{$#1$}}}\raisebox{-.6ex}{$\begin{array}[b]{*{1}c}\cline{1-1}\lr{\phantom{x}}\\\cline{1-1}\end{array}$}}$};
\node (node_2) at (127.0bp,121.5bp) [draw,draw=none] {${\def\lr#1{\multicolumn{1}{|@{\hspace{.6ex}}c@{\hspace{.6ex}}|}{\raisebox{-.3ex}{$#1$}}}\raisebox{-.6ex}{$\begin{array}[b]{*{1}c}\cline{1-1}\lr{\phantom{x}}\\\cline{1-1}\lr{\phantom{x}}\\\cline{1-1}\end{array}$}}$};
\node (node_6) at (169.0bp,121.5bp) [draw,draw=none] {${\def\lr#1{\multicolumn{1}{|@{\hspace{.6ex}}c@{\hspace{.6ex}}|}{\raisebox{-.3ex}{$#1$}}}\raisebox{-.6ex}{$\begin{array}[b]{*{2}c}\cline{1-2}\lr{\phantom{x}}&\lr{\phantom{x}}\\\cline{1-2}\end{array}$}}$};
\node (node_3) at (99.0bp,194.5bp) [draw,draw=none] {${\def\lr#1{\multicolumn{1}{|@{\hspace{.6ex}}c@{\hspace{.6ex}}|}{\raisebox{-.3ex}{$#1$}}}\raisebox{-.6ex}{$\begin{array}[b]{*{1}c}\cline{1-1}\lr{\phantom{x}}\\\cline{1-1}\lr{\phantom{x}}\\\cline{1-1}\lr{\phantom{x}}\\\cline{1-1}\end{array}$}}$};
\node (node_7) at (148.0bp,194.5bp) [draw,draw=none] {${\def\lr#1{\multicolumn{1}{|@{\hspace{.6ex}}c@{\hspace{.6ex}}|}{\raisebox{-.3ex}{$#1$}}}\raisebox{-.6ex}{$\begin{array}[b]{*{2}c}\cline{1-2}\lr{\phantom{x}}&\lr{\phantom{x}}\\\cline{1-2}\lr{\phantom{x}}\\\cline{1-1}\end{array}$}}$};
\node (node_4) at (51.0bp,279.0bp) [draw,draw=none] {${\def\lr#1{\multicolumn{1}{|@{\hspace{.6ex}}c@{\hspace{.6ex}}|}{\raisebox{-.3ex}{$#1$}}}\raisebox{-.6ex}{$\begin{array}[b]{*{1}c}\cline{1-1}\lr{\phantom{x}}\\\cline{1-1}\lr{\phantom{x}}\\\cline{1-1}\lr{\phantom{x}}\\\cline{1-1}\lr{\phantom{x}}\\\cline{1-1}\end{array}$}}$};
\node (node_8) at (99.0bp,279.0bp) [draw,draw=none] {${\def\lr#1{\multicolumn{1}{|@{\hspace{.6ex}}c@{\hspace{.6ex}}|}{\raisebox{-.3ex}{$#1$}}}\raisebox{-.6ex}{$\begin{array}[b]{*{2}c}\cline{1-2}\lr{\phantom{x}}&\lr{\phantom{x}}\\\cline{1-2}\lr{\phantom{x}}\\\cline{1-1}\lr{\phantom{x}}\\\cline{1-1}\end{array}$}}$};
\node (node_5) at (9.0bp,375.0bp) [draw,draw=none] {${\def\lr#1{\multicolumn{1}{|@{\hspace{.6ex}}c@{\hspace{.6ex}}|}{\raisebox{-.3ex}{$#1$}}}\raisebox{-.6ex}{$\begin{array}[b]{*{1}c}\cline{1-1}\lr{\phantom{x}}\\\cline{1-1}\lr{\phantom{x}}\\\cline{1-1}\lr{\phantom{x}}\\\cline{1-1}\lr{\phantom{x}}\\\cline{1-1}\lr{\phantom{x}}\\\cline{1-1}\end{array}$}}$};
\node (node_9) at (51.0bp,375.0bp) [draw,draw=none] {${\def\lr#1{\multicolumn{1}{|@{\hspace{.6ex}}c@{\hspace{.6ex}}|}{\raisebox{-.3ex}{$#1$}}}\raisebox{-.6ex}{$\begin{array}[b]{*{2}c}\cline{1-2}\lr{\phantom{x}}&\lr{\phantom{x}}\\\cline{1-2}\lr{\phantom{x}}\\\cline{1-1}\lr{\phantom{x}}\\\cline{1-1}\lr{\phantom{x}}\\\cline{1-1}\end{array}$}}$};
\node (node_12) at (204.0bp,194.5bp) [draw,draw=none] {${\def\lr#1{\multicolumn{1}{|@{\hspace{.6ex}}c@{\hspace{.6ex}}|}{\raisebox{-.3ex}{$#1$}}}\raisebox{-.6ex}{$\begin{array}[b]{*{3}c}\cline{1-3}\lr{\phantom{x}}&\lr{\phantom{x}}&\lr{\phantom{x}}\\\cline{1-3}\end{array}$}}$};
\node (node_10) at (148.0bp,279.0bp) [draw,draw=none] {${\def\lr#1{\multicolumn{1}{|@{\hspace{.6ex}}c@{\hspace{.6ex}}|}{\raisebox{-.3ex}{$#1$}}}\raisebox{-.6ex}{$\begin{array}[b]{*{2}c}\cline{1-2}\lr{\phantom{x}}&\lr{\phantom{x}}\\\cline{1-2}\lr{\phantom{x}}&\lr{\phantom{x}}\\\cline{1-2}\end{array}$}}$};
\node (node_13) at (204.0bp,279.0bp) [draw,draw=none] {${\def\lr#1{\multicolumn{1}{|@{\hspace{.6ex}}c@{\hspace{.6ex}}|}{\raisebox{-.3ex}{$#1$}}}\raisebox{-.6ex}{$\begin{array}[b]{*{3}c}\cline{1-3}\lr{\phantom{x}}&\lr{\phantom{x}}&\lr{\phantom{x}}\\\cline{1-3}\lr{\phantom{x}}\\\cline{1-1}\end{array}$}}$};
\node (node_11) at (98.0bp,375.0bp) [draw,draw=none] {${\def\lr#1{\multicolumn{1}{|@{\hspace{.6ex}}c@{\hspace{.6ex}}|}{\raisebox{-.3ex}{$#1$}}}\raisebox{-.6ex}{$\begin{array}[b]{*{2}c}\cline{1-2}\lr{\phantom{x}}&\lr{\phantom{x}}\\\cline{1-2}\lr{\phantom{x}}&\lr{\phantom{x}}\\\cline{1-2}\lr{\phantom{x}}\\\cline{1-1}\end{array}$}}$};
\node (node_14) at (151.0bp,375.0bp) [draw,draw=none] {${\def\lr#1{\multicolumn{1}{|@{\hspace{.6ex}}c@{\hspace{.6ex}}|}{\raisebox{-.3ex}{$#1$}}}\raisebox{-.6ex}{$\begin{array}[b]{*{3}c}\cline{1-3}\lr{\phantom{x}}&\lr{\phantom{x}}&\lr{\phantom{x}}\\\cline{1-3}\lr{\phantom{x}}\\\cline{1-1}\lr{\phantom{x}}\\\cline{1-1}\end{array}$}}$};
\node (node_15) at (209.0bp,375.0bp) [draw,draw=none] {${\def\lr#1{\multicolumn{1}{|@{\hspace{.6ex}}c@{\hspace{.6ex}}|}{\raisebox{-.3ex}{$#1$}}}\raisebox{-.6ex}{$\begin{array}[b]{*{3}c}\cline{1-3}\lr{\phantom{x}}&\lr{\phantom{x}}&\lr{\phantom{x}}\\\cline{1-3}\lr{\phantom{x}}&\lr{\phantom{x}}\\\cline{1-2}\end{array}$}}$};
\node (node_16) at (273.0bp,279.0bp) [draw,draw=none] {${\def\lr#1{\multicolumn{1}{|@{\hspace{.6ex}}c@{\hspace{.6ex}}|}{\raisebox{-.3ex}{$#1$}}}\raisebox{-.6ex}{$\begin{array}[b]{*{4}c}\cline{1-4}\lr{\phantom{x}}&\lr{\phantom{x}}&\lr{\phantom{x}}&\lr{\phantom{x}}\\\cline{1-4}\end{array}$}}$};
\node (node_17) at (273.0bp,375.0bp) [draw,draw=none] {${\def\lr#1{\multicolumn{1}{|@{\hspace{.6ex}}c@{\hspace{.6ex}}|}{\raisebox{-.3ex}{$#1$}}}\raisebox{-.6ex}{$\begin{array}[b]{*{4}c}\cline{1-4}\lr{\phantom{x}}&\lr{\phantom{x}}&\lr{\phantom{x}}&\lr{\phantom{x}}\\\cline{1-4}\lr{\phantom{x}}\\\cline{1-1}\end{array}$}}$};
\node (node_18) at (348.0bp,375.0bp) [draw,draw=none] {${\def\lr#1{\multicolumn{1}{|@{\hspace{.6ex}}c@{\hspace{.6ex}}|}{\raisebox{-.3ex}{$#1$}}}\raisebox{-.6ex}{$\begin{array}[b]{*{5}c}\cline{1-5}\lr{\phantom{x}}&\lr{\phantom{x}}&\lr{\phantom{x}}&\lr{\phantom{x}}&\lr{\phantom{x}}\\\cline{1-5}\end{array}$}}$};
\draw [black,<->] (node_0) ..controls (148.0bp,21.694bp) and (148.0bp,31.939bp)  .. (node_1);
\draw [black,<->] (node_1) ..controls (142.45bp,77.079bp) and (138.92bp,87.018bp)  .. (node_2);
\draw [black,<->] (node_1) ..controls (154.08bp,78.572bp) and (158.67bp,91.478bp)  .. (node_6);
\draw [black,<->] (node_2) ..controls (118.21bp,144.79bp) and (114.45bp,154.31bp)  .. (node_3);
\draw [black,<->] (node_2) ..controls (134.08bp,146.45bp) and (137.67bp,158.58bp)  .. (node_7);
\draw [black,<->] (node_3) ..controls (83.043bp,222.93bp) and (73.343bp,239.6bp)  .. (node_4);
\draw [black,<->] (node_3) ..controls (99.0bp,225.67bp) and (99.0bp,236.9bp)  .. (node_8);
\draw [black,<->] (node_4) ..controls (36.163bp,313.21bp) and (28.682bp,329.95bp)  .. (node_5);
\draw [black,<->] (node_4) ..controls (51.0bp,315.96bp) and (51.0bp,327.25bp)  .. (node_9);
\draw [black,<->] (node_6) ..controls (163.66bp,140.56bp) and (159.14bp,155.83bp)  .. (node_7);
\draw [black,<->] (node_6) ..controls (178.7bp,142.17bp) and (188.12bp,161.29bp)  .. (node_12);
\draw [black,<->] (node_7) ..controls (132.89bp,220.94bp) and (124.09bp,235.76bp)  .. (node_8);
\draw [black,<->] (node_7) ..controls (148.0bp,222.1bp) and (148.0bp,239.01bp)  .. (node_10);
\draw [black,<->] (node_7) ..controls (166.45bp,222.68bp) and (178.68bp,240.69bp)  .. (node_13);
\draw [black,<->] (node_8) ..controls (82.658bp,312.0bp) and (75.501bp,326.02bp)  .. (node_9);
\draw [black,<->] (node_8) ..controls (98.646bp,313.27bp) and (98.476bp,329.25bp)  .. (node_11);
\draw [black,<->] (node_8) ..controls (117.61bp,313.64bp) and (126.77bp,330.2bp)  .. (node_14);
\draw [black,<->] (node_10) ..controls (133.16bp,307.9bp) and (122.63bp,327.7bp)  .. (node_11);
\draw [black,<->] (node_10) ..controls (167.29bp,309.73bp) and (182.6bp,333.32bp)  .. (node_15);
\draw [black,<->] (node_12) ..controls (204.0bp,215.5bp) and (204.0bp,236.43bp)  .. (node_13);
\draw [black,<->] (node_12) ..controls (222.38bp,217.48bp) and (244.5bp,243.92bp)  .. (node_16);
\draw [black,<->] (node_13) ..controls (188.2bp,308.02bp) and (176.89bp,328.07bp)  .. (node_14);
\draw [black,<->] (node_13) ..controls (205.56bp,309.31bp) and (206.77bp,332.03bp)  .. (node_15);
\draw [black,<->] (node_13) ..controls (225.92bp,309.86bp) and (243.46bp,333.75bp)  .. (node_17);
\draw [black,<->] (node_16) ..controls (273.0bp,302.25bp) and (273.0bp,329.02bp)  .. (node_17);
\draw [black,<->] (node_16) ..controls (292.37bp,304.27bp) and (318.72bp,337.3bp)  .. (node_18);
\end{tikzpicture}$}
.
\end{gather*}
We get $\pfdim\Gamma(Y)=\infty$ since $\Gamma(Y)$ contains star graphs 
for arbitrarily large $N$. To see this take the vertex $0$ of the star 
graph to be a staircase partition $(t,t-1,t-2,\dots,1)$ for some large enough $t$.
\end{Example}

\begin{Proposition}\label{P:GrowthRateSuper}
The following categories are of superexponential growth.
\begin{enumerate}

\item The category $\rep\big(\mathrm{SL}_{\N}(\C),\C^{\N}\big)$.

\item The Deligne category $\rep(S_{t},\C)$, where $t\in\C\setminus\Z$.

\item The so-called \emph{Delannoy categories} from \cite{HarSnoSny-delannoycat}.

\end{enumerate}
\end{Proposition}

\begin{proof}
We can apply \autoref{E:GrowthRateStarGraph} three times:
For $\rep\big(\mathrm{SL}_{\N}(\C),\C^{\N}\big)$ this follows directly from \autoref{T:GrowthRateMainTheorem}.(a) and \autoref{E:GrowthRateStarGraph}. 

The following argument is based on the description of the fusion graph of $\rep(S_{t},\C)$ for $t\in\C\setminus\Z$, which was provided to us by Victor Ostrik via email. We thank Victor for providing us with this description.
For $\rep(S_{t},\C)$ we recall that the combinatorics of tensoring in
is very similar to the combinatorics of tensoring with $\star_{N}$-representations 
if $N\gg 0$. The correspondence is as follows: if you have a partition
$(p_{1},p_{2},\dots)$ labeling object in the Deligne category, it corresponds
to a partition $(N-(p_{1}+p_{2}+\dots),p_{1},p_{2},\dots)$ for $\star_{N}$.
Now tensoring
by a 1-node Young diagram object is given by removing one node whenever possible and add it back. This description implies the existence of star 
graphs, so we are done in the same way as for $\rep\big(\mathrm{SL}_{\N}(\C),\C^{\N}\big)$.

For the Delannoy categories from \cite{HarSnoSny-delannoycat} we can again use star graphs for the fusion graphs computed in \cite[Section 7]{HarSnoSny-delannoycat}.
\end{proof}

\begin{Remark}\label{R:GrowthRateKKO}
By \cite{GrTu}, ``most'' of the categories discussed in \cite{KKO-twodimTQFT} exhibit superexponential growth, and we suspect that the same is true for the generalizations of the Delannoy categories from \cite{HarSnoSny-delannoycat}.
\end{Remark}

%
%
%
%
%
%
%


\section{Recurrent and transient growth problems}\label{S:RT}


We retain the setup in \autoref{S:Basics}.
Following the 
classical theory of random walks, see {\eg} 
\cite{Ki-symbolic-dynamics} where the reader will also find some standard terminology,
we now define \emph{recurrent and transient growth problems}.

To this end, recall that a vertex $v\in\Gamma$ in a graph with countably many vertices
(not just a fusion graph) is \emph{recurrent} if the probability 
of returning to $v$ in the random walk with all edge labels one on $\Gamma$ 
is $1$, and \emph{transient} if its not recurrent.
If $v\in\Gamma$ is recurrent, then we call it \emph{positive recurrent}
if the expected amount of time between recurrences is finite, and 
\emph{null recurrent} otherwise.

\begin{Remark}
Technically speaking, we are slightly abusing the terminology here. What we refer to as a `random walk' does not align with the standard definition, since we often assign a weight of one to each edge. In a true random walk, each edge would be given a label such that the sum of all edge weights equals one.
\end{Remark}

\begin{Lemma}\label{L:RTConnectedComponents}
For a fusion graph $\Gamma$ we have:
\begin{enumerate}

\item Every vertex $v$ of $\Gamma$ is either recurrent or transient.

\item Every recurrent vertex of $\Gamma$ is either positively recurrent or null recurrent.

\item Let $C(\Gamma)$ be a connected component. If $v\in C(\Gamma)$ is positively recurrent (or null recurrent or transient), then so is any other $w\in C(\Gamma)$.

\end{enumerate}

\end{Lemma}

\begin{proof}
(a) and (b) are true by definition (and only stressed because of the many alternative definitions of recurrent and transient, which are not immediately opposite of each other), and (c) is classical.
\end{proof}

Note that \autoref{L:RTConnectedComponents} allows us to consider (positively/null) recurrence and transience for 
one fixed connected component. We will use this below.

To state the main definition, we recall the notion of 
PF dimension 
of irreducible (not necessary finite) $\Rplus$-matrices {\`a} la Vere--Jones \cite{VeJo-ergodic-properties-one} (justifying the subscript in the following notation). For such a matrix $M=(m_{ij})_{i,j\in I}$ 
for a countable set $I$, let
\begin{gather}\label{Eq:RecurrentPFDim}
\pfdimtwo M=\lim_{n\to\infty}\sqrt[h\cdot n]{m_{ij}^{(h\cdot n)}},
\end{gather}
where $m_{ij}^{(n)}$ is the $(i,j)$-entry of $M^{n}$, and $h\in\Z_{\geq 1}$ 
is the period. By \cite[Theorem A]{VeJo-ergodic-properties-one} $\pfdimtwo M$ is independent of $i$ and $j$.

Recalling that irreducible matrices correspond to strongly connected graphs, we can define:

\begin{Definition}\label{D:RecurrentPFDim}
For a strongly connected graph $\Gamma$ with countably many vertices define
the \emph{PF dimension}
$\pfdimtwo\Gamma$ using \autoref{Eq:RecurrentPFDim} with $M$ being 
the adjacency matrix of $\Gamma$.
\end{Definition}

\begin{Remark}
A generalization of $\pfdimtwo\Gamma$ to graphs that are not necessarily strongly connected was worked out in \cite{Tw-nonnegative-matrices}. It however turns out that this generalization is not useful for our purposes.
\end{Remark}

The following is immediate:

\begin{Lemma}\label{L:RecurrentPFFinite}
If $|I|<\infty$, then 
\autoref{D:RecurrentPFDim} agrees with the usual 
definition of PF dimension using the largest eigenvalue of the associated matrix (or, equivalently, graph).\qed
\end{Lemma}

\begin{Remark}
We use the notation $\pfdimtwo$ to indicate that this PF dimension is different from the 
one we introduced in \autoref{S:GrowthRate}. They however agree 
if $|I|<\infty$ by \autoref{L:GrowthRatePFFinite} and \autoref{L:RecurrentPFFinite}.
\end{Remark}

\begin{Example}\label{E:RecurrentPFDim}
We continue with \autoref{E:GrowthRateSLTwoOne}. In this case we get
$\pfdimtwo\Gamma(\mathrm{SL}_{2})=2$ since
\begin{gather*}
m_{11}^{(2n)}(\mathrm{SL}_{2})
=\frac{1}{n+1}\binom{2n}{n}\sim\pi^{-1/2}\cdot n^{-3/2}\cdot 2^{2n}.
\end{gather*}
This equation follows from a standard argument. Note that the period is two, 
which is why we consider every second value of $m_{11}^{(n)}(\mathrm{SL}_{2})$ only.
\end{Example}

We also recall the notion of
\emph{final basic classes (FBC)}. Firstly, a \emph{class} $C(\Gamma)$ is a strongly connected component of $\Gamma$ and such a class is \emph{basic} if:
\begin{gather*}
\pfdimtwo C(\Gamma)\geq\pfdimtwo C^{\prime}(\Gamma)
\text{ for all strongly connected components }C^{\prime}(\Gamma).
\end{gather*}
Finally, a basic class is \emph{final} if there is no path to any other basic class.
(Note that a FBC does not need to be final in $\Gamma$ itself.)

\begin{Example}\label{E:RecurrentPFDimEasy}
For a strongly connected graph $\Gamma$ is the only class, so it is a FBC. In particular, 
what we will see below is a generalization of the strongly connected situation.
\end{Example}

\begin{Example}\label{E:RecurrentSL2Fp}
Consider $G=\mathrm{SL}_{2}(\bF[2])$ and take the defining $G$-representation 
$\obstuff{X}=\bF[2]^{2}$ over the field $\bF[2]$. The growth problem
$(G,\obstuff{X})$ has the associated fusion graph (we illustrate a cutoff):
\begin{gather*}
\Gamma\big(\mathrm{SL}_{2}(\bF[2])\big)
=
\begin{tikzpicture}[anchorbase]
\node at (0,0) {\includegraphics[height=6.0cm]{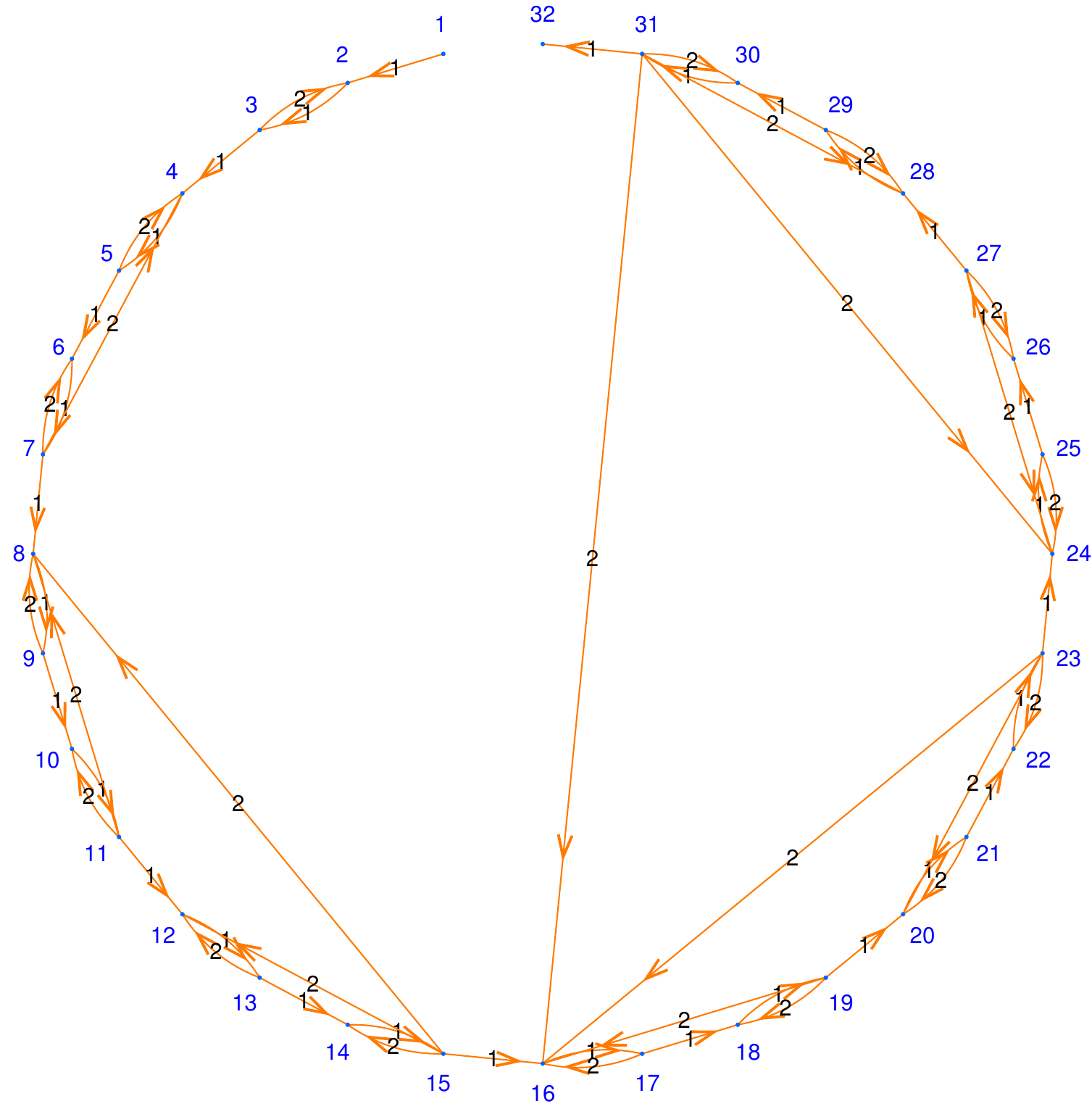}};
\end{tikzpicture}
.
\end{gather*}
The pattern is as follows: The vertices are labeled with $\Z_{\geq 1}$. There is always one forward arrow, from vertex $i$ to vertex $i+1$. The backward arrows 
are of weight $2$ and start at $2^{k}-1$ and end at $2^{k-2}$ for all $k\in\Z_{\geq 2}$.
See, for example, \cite{SuTuWeZh-mixed-tilting} for details.

The strongly connected components are indexed by $2^{k}$ for $k\in\N$, 
with $2^{k}$ being their minimal vertex.
In this case there is no FBC since all strongly connected components 
have PF dimension $<2$, but the limit $k\to\infty$ of their PF dimensions is approaching $2$. 

In a bit more details, one can check (by calculation) that the PF dimension of the cutoffs will be
\begin{gather*}
0,\sqrt{2},\sqrt{2+\sqrt{2}},\sqrt{2+\sqrt{2+\sqrt{2}}}\text{ {\etc}}
\end{gather*}
where the jumps happen when the number of vertices increases from $2^{k}$ to $2^{k}+1$ in the sequence of naive cutoffs. The limit approaches $2$, by classical results about sequences of nested radicals.
\end{Example}

The following is a key definition:

\begin{Definition}\label{D:RTMain}
Let $(R,c)$ be a growth problem, and let FBC refer to its fusion graph $\Gamma$.
\begin{enumerate}

\item The growth problem $(R,c)$ is \emph{recurrent} if all of its 
FBCs are recurrent and there exists at least one FBC. The 
growth problem is \emph{transient} otherwise.

\item If $(R,c)$ is recurrent, then we say $(R,c)$ is \emph{positively recurrent} if 
all of its FBCs are positively recurrent,  The 
growth problem is \emph{null recurrent} otherwise.

\end{enumerate}
We use the same terminology for graphs in general.
\end{Definition}

Note that a transient FBC dominates recurrent ones in the sense that the presence of even a single transient FBC renders the growth problem transient, regardless of how many FBCs are recurrent. Similarly, a null recurrent FBC dominates positively recurrent ones.

\begin{Lemma}\label{L:RTMain}
Every growth problem is either recurrent or transient.
Moreover, every recurrent growth problem is either positively or null recurrent.
\end{Lemma}

\begin{proof}
Immediate by definition (but stated for the same reasons as \autoref{L:RTConnectedComponents}).
\end{proof}

\begin{Example}\label{E:GrowthRateSLTwoThree}
Take the graph 
\begin{gather*}
\Gamma(\Z)
=
\begin{tikzpicture}[anchorbase]
\begin{scope}[every node/.style={circle,inner sep=0pt,text width=6mm,align=center,draw=white,fill=white}]
\node (B) at (2,0) {$\bullet$};
\node (C) at (4,0) {$\bullet$};
\node (D) at (6,0) {$\bullet$};
\node (E) at (8,0) {$\bullet$};
\node (F) at (10,0) {$\bullet$};
\end{scope}
\node (A) at (0,0) {$\dots$};
\node (G) at (12,0) {$\dots$};
\begin{scope}[>={Stealth[black]},
every edge/.style={draw=orchid,very thick}]
\path [->] ($(A)+(0.25,0.075)$) edge ($(B)+(-0.25,0.075)$);
\path [->] ($(B)+(0.25,0.075)$) edge ($(C)+(-0.25,0.075)$);
\path [->] ($(C)+(0.25,0.075)$) edge ($(D)+(-0.25,0.075)$);
\path [->] ($(D)+(0.25,0.075)$) edge ($(E)+(-0.25,0.075)$);
\path [->] ($(E)+(0.25,0.075)$) edge ($(F)+(-0.25,0.075)$);
\path [->] ($(F)+(0.25,0.075)$) edge ($(G)+(-0.25,0.075)$);
\path [<-] ($(A)+(0.25,-0.075)$) edge ($(B)+(-0.25,-0.075)$);
\path [<-] ($(B)+(0.25,-0.075)$) edge ($(C)+(-0.25,-0.075)$);
\path [<-] ($(C)+(0.25,-0.075)$) edge ($(D)+(-0.25,-0.075)$);
\path [<-] ($(D)+(0.25,-0.075)$) edge ($(E)+(-0.25,-0.075)$);
\path [<-] ($(E)+(0.25,-0.075)$) edge ($(F)+(-0.25,-0.075)$);
\path [<-] ($(F)+(0.25,-0.075)$) edge ($(G)+(-0.25,-0.075)$);
\end{scope}
\end{tikzpicture}
.
\end{gather*}
The associated random walk is the classical random walk on $\Z$. This is Polya's
first example of a recurrent graph 
\cite{Po-random-walks}. Additionally, this example is null recurrent.

In contrast, the random walk on $\Gamma(\mathrm{SL}_{2})$
from \autoref{Eq:GrowthRateSLTwoOne}
is transient. Below we will see that this 
is no coincidence, {\cf} \autoref{P:RTReductive}.

The difference becomes evident when one looks at the number of path of length $n$ starting at the origin, and ending at vertex $v$. For $n=200$, plotting this in $\R^2$ with $(\text{end vertex},\text{number of paths})$ gives:
\begin{gather*}
\Z\colon
\begin{tikzpicture}[anchorbase,scale=1]
\node at (0,0) {\includegraphics[height=3.6cm]{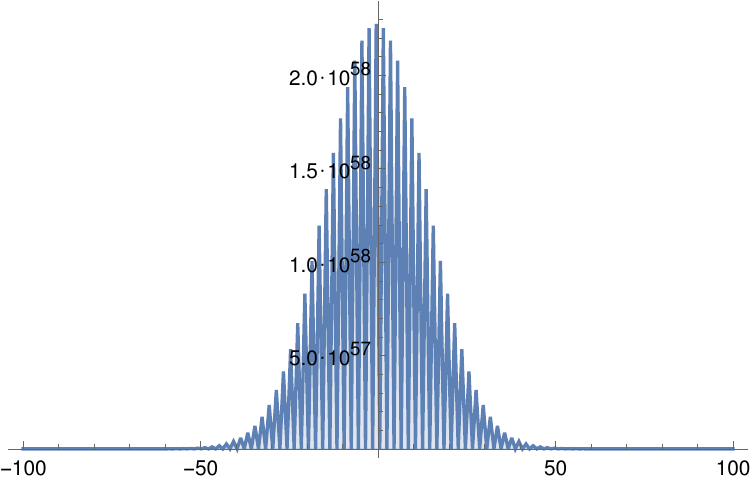}};
\end{tikzpicture}
,\quad
\N\colon
\begin{tikzpicture}[anchorbase,scale=1]
\node at (0,0) {\includegraphics[height=3.6cm]{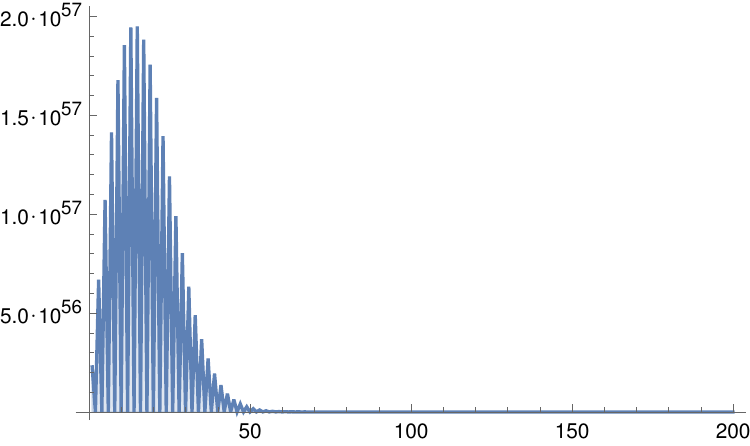}};
\end{tikzpicture}
.
\end{gather*}
Note that for $\Gamma(\mathrm{SL}_{2})\cong\N$ the peak of the binomial distribution is roughly at $\sqrt{200}\approx 14.1$. In fact, the endpoints of paths in this case move to infinity, while 
they stay at the origin for $\Z$.

For completeness:
The graph $\Gamma(\Z)$ is also the fusion graph for a monoidal category: take $(\rep(\mathbb{S}^{1}),\obstuff{X}_{1}\oplus\obstuff{X}_{-1})$ with $\obstuff{X}_{\pm 1}$ being the rotation of an irrational angle $\pm\theta$.
\end{Example}

\begin{Example}\label{E:GrowthRateSLThree}
Similarly as for $\mathrm{SL}_{2}(\C)$, the growth problem for $G=\mathrm{SL}_{3}(\C)$ and $\obstuff{X}=\C^{3}$ with the 
defining $G$-action is also transient. In this case, the fusion graph is (here a cutoff):
\begin{gather*}
\Gamma(\mathrm{SL}_{3})\leftrightsquigarrow
\begin{tikzpicture}[anchorbase,scale=1]
\node at (0,0) {\includegraphics[height=3.6cm]{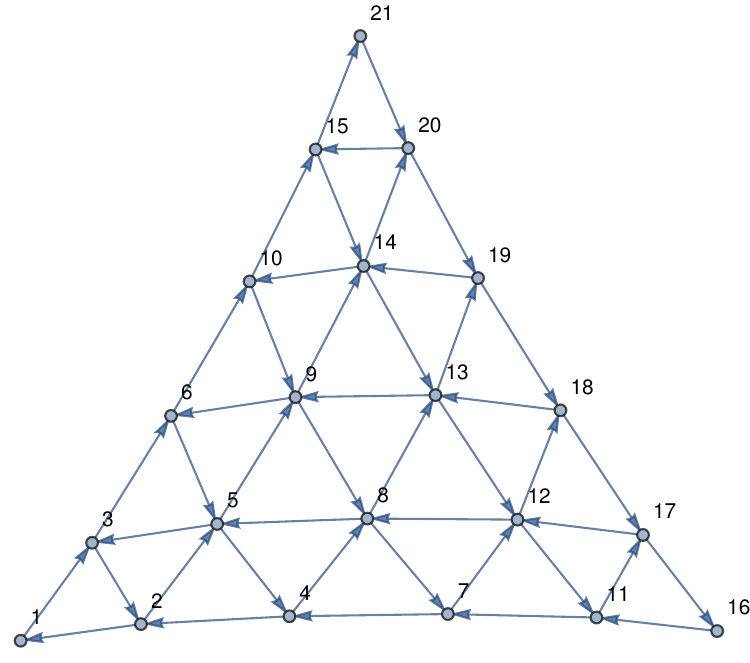}};
\end{tikzpicture}
.
\end{gather*}
And similarly as in \autoref{E:GrowthRateSLTwoThree}, the number of paths 
of length $n$ starting at $1$ and ending at $v$ moves outwards (illustrated using a large enough cutoff):
\begin{gather*}
n=50\colon
\begin{tikzpicture}[anchorbase,scale=1]
\node at (0,0) {\includegraphics[height=2.0cm]{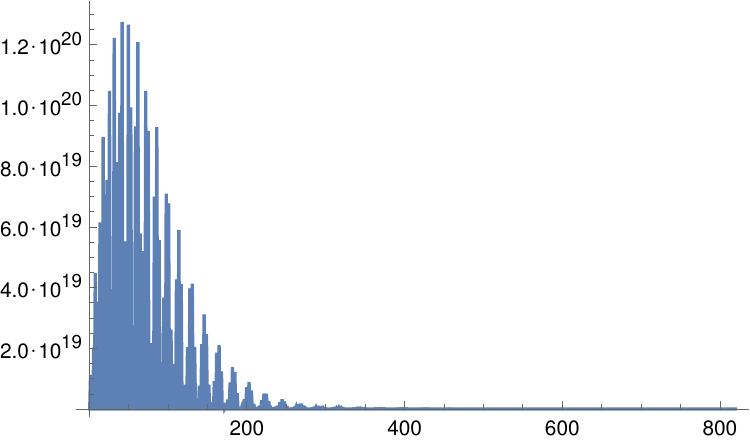}};
\end{tikzpicture}
,\quad
n=100\colon
\begin{tikzpicture}[anchorbase,scale=1]
\node at (0,0) {\includegraphics[height=2.0cm]{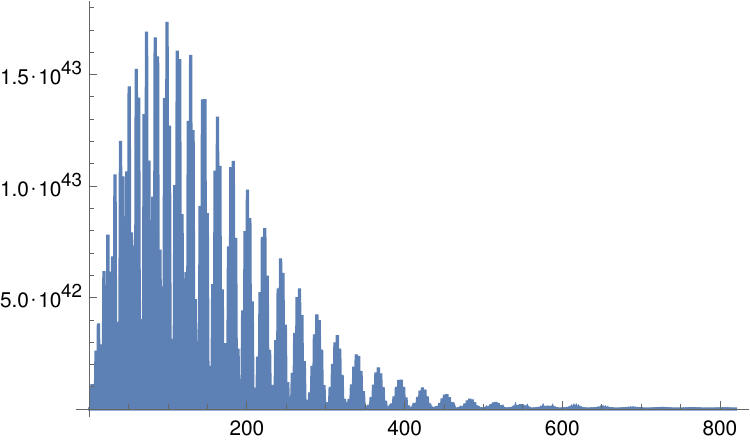}};
\end{tikzpicture}
,\quad
n=200\colon
\begin{tikzpicture}[anchorbase,scale=1]
\node at (0,0) {\includegraphics[height=2.0cm]{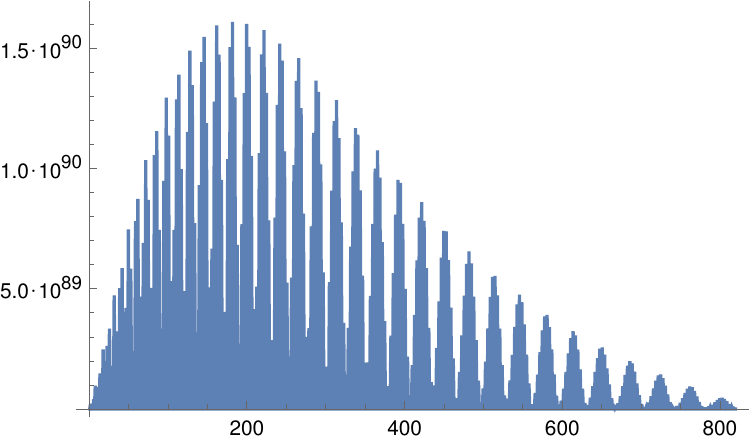}};
\end{tikzpicture}
.
\end{gather*}
In \autoref{P:RTReductive} we will see that this is not a coincidence.
\end{Example}

\begin{Example}\label{E:GrowthRateFinite}
Every finite growth problem is positively recurrent. More generally, every growth problem 
where the FBCs are finite is also positively recurrent. Explicitly, let $G=\Z/2\Z\times\Z/2\Z$ be the \emph{Klein four group}, and we take the ground field $\K=\F[2]$. Then one has $\K[G]\cong\K[X,Y]/(X^{2},Y^{2})$, and using this presentation one can check that
\begin{gather*}
\inde[3]\leftrightsquigarrow
\begin{tikzpicture}[anchorbase]
\begin{scope}[every node/.style={circle,inner sep=0pt,text width=6mm,align=center,draw=white,fill=white}]
\node (A) at (0,0) {$\bullet$};
\node (B) at (-1,-1) {$\bullet$};
\node (C) at (1,-1) {$\bullet$};
\end{scope}
\begin{scope}[>={Stealth[black]},
every edge/.style={draw=black,thick}]
\path [->] (A) edge node[left,xshift=-0.15cm]{$X$} (B);
\path [->] (A) edge node[right,xshift=0.15cm]{$Y$} (C);
\end{scope}
\end{tikzpicture}
\end{gather*}	
defines an indecomposable $G$-representation of dimension three, call it $V=\inde[3]$. Moreover, it is also easy to see that 
$G$ has indecomposable $G$-representations $\inde[k]$, for odd
$k\in\Z_{\geq 1}$, where $\inde[k]$ is the unique indecomposable summand of $V^{\otimes(k-1)/2}$ that is not projective.
Note that $\inde[1]\cong\simple[1]$ is the trivial $G$-representation. One can check that $\dim_{\K}\inde[k]=k$ and that all odd numbers appear.

Let $\prin[4]$ denote the regular $G$-representation.
An easy calculation gives us then the following fusion graph.
\begin{gather}\label{Eq:GrowthRateMainTheoremKlein}
\begin{gathered}
\Gamma(Kl)=
\begin{tikzpicture}[anchorbase]
\begin{scope}[every node/.style={circle,inner sep=0pt,text width=6mm,align=center,draw=black,fill=white}]
\node (A) at (0,0) {$\simple[1]$};
\node (B) at (2,0) {$\inde[3]$};
\node (C) at (4,0) {$\inde[5]$};
\node (D) at (6,0) {$\inde[7]$};
\node (E) at (8,0) {$\inde[9]$};
\node (F) at (10,0) {$\inde[11]$};
\node (P) at (6,-2) {$\prin[4]$};
\end{scope}
\node (G) at (12,0) {$\dots$};
\begin{scope}[>={Stealth[black]},
every edge/.style={draw=orchid,very thick}]
\path [->] (A) edge (B);
\path [->] (B) edge (C);
\path [->] (C) edge (D);
\path [->] (D) edge (E);
\path [->] (E) edge (F);
\path [->] (F) edge (G); 
\path [->] (B) edge node[left,xshift=-0.15cm]{$1$} (P);
\path [->] (C) edge node[left,xshift=-0.15cm]{$2$} (P);
\path [->] (D) edge node[left,xshift=-0.15cm]{$3$} (P);
\path [->] (E) edge node[left,xshift=-0.15cm]{$4$} (P);
\path [->] (F) edge node[left,xshift=-0.15cm]{$5$} (P);
\path [->] (G) edge node[left,xshift=-0.15cm]{$\dots$} (P);
\path [->] (P) edge [loop below] node{$3$} (P);
\end{scope}
\end{tikzpicture}
,
\\
M(Kl)=
\scalebox{0.5}{$\left(
\begin{array}{cccccccc}
\cellcolor{spinach!25}3 & 0 & 1 & 2 & 3 & 4 & 5 & \dots \\
0 & 0 & 0 & 0 & 0 & 0 & 0 & \dots \\
0 & 1 & 0 & 0 & 0 & 0 & 0 & \dots \\
0 & 0 & 1 & 0 & 0 & 0 & 0 & \dots \\
0 & 0 & 0 & 1 & 0 & 0 & 0 & \dots \\
0 & 0 & 0 & 0 & 1 & 0 & 0 & \dots \\
0 & 0 & 0 & 0 & 0 & 1 & 0 & \dots \\
\vdots & \vdots & \vdots & \vdots & \vdots & \vdots & \vdots & \ddots \\
\end{array}
\right)$}
.
\end{gathered}
\end{gather}
The labeled loop at the vertex with label $3$ represents three loops.

Let us consider the growth rate problem for $(G,V)$. Then this is positively recurrent 
since the only FBC is $\{\prin[4]\}$ and therefore a finite graph.
\end{Example}

For a more general statement see \autoref{P:RTReductive} below.

\begin{Example}\label{E:RTSLTwo}
The growth problem in \autoref{E:RecurrentSL2Fp} is transient since 
there is no FBC.
\end{Example}

Recall that a group $G$ is \emph{virtually} $H$ if 
$H$ is isomorphic to a finite index subgroup of $G$.

\begin{Theorem}[Polya's random walk theorem for representations]\label{T:RTPolya}
Let $\K$ be a field of characteristic zero.
For any (finitely presentable) group $G$ and finite dimensional faithful 
completely reducible $G$-representation $\obstuff{X}$ the following are equivalent:
\begin{enumerate}

\item The growth problem $(G,\obstuff{X})$ is recurrent.

\item The Zariski closure of the image of $G$ in $\mathrm{GL}(\obstuff{X})$ is virtually a torus of rank $0$, $1$ or $2$. \label{I:RTPolya}

\end{enumerate}
Moreover, the following are also equivalent:
\begin{enumerate}[label=(\roman*)]

\item The growth problem $(G,\obstuff{X})$ is positively recurrent.

\item The Zariski closure of the image of $G$ in $\mathrm{GL}(\obstuff{X})$ is virtually a torus of rank $0$, i.e. $G$ is finite.

\end{enumerate}
If $\obstuff{X}$ is not faithful, then go to the biggest 
quotient group that acts faithfully and repeat.
\end{Theorem}

\begin{proof}
We will make use of the following. We will formulated it for bialgebra, 
but we use the same notation as for group.

\begin{Lemma}\label{L:RTReductiveOne}
Let $G$ be a bialgebra, and $H\subset G$ be a subbialgebra. Let $\obstuff{X}$
be a $G$-representation such that:
\begin{enumerate}

\item The growth problem $(G,\obstuff{X})$ has a strongly connected fusion graph with $\pfdimtwo\obstuff{X}<\infty$.

\item The growth problem $(H,\obstuff{X})$ obtained via restriction has a strongly connected fusion graph with the same PF dimension $\pfdimtwo\obstuff{X}$.
    
\end{enumerate}
Let $\obstuff{X}\cong\bigoplus_{i=1}^{j}\obstuff{X}_{i}$ as an $H$-representation obtained by restriction.
If the growth problem $(H,\obstuff{X}_{i})$ is transient for some $i\in\{1,\dots,j\}$, then $(G,\obstuff{X})$ is also transient.
\end{Lemma}

\begin{proof}
We only prove the statement in the aperiodic case, to keep the formulas simple.
First note that the transience of $(H,\obstuff{X}_{i})$ implies that $(H,\obstuff{X})$ is also transient as there is a subgraph corresponding to $(H,\obstuff{X}_{i})$. Now, since $(H,\obstuff{X})$ is transient with strongly connected fusion graph we have
\begin{gather*}
\sum_{n\in\N}h_{11}^{(n)}\cdot (\pfdimtwo\obstuff{X})^{-n}<\infty,
\end{gather*}
where $h_{11}^{(n)}$ is the $(1,1)$-entry 
of the $n$th power of the action matrix corresponding to the fusion graph of $(H,\obstuff{X})$. Using a similar notation for $(G,\obstuff{X})$, we also have
\begin{gather*}
\sum_{n\in\N}g_{11}^{(n)}\cdot (\pfdimtwo\obstuff{X})^{-n}
\leq
\sum_{n\in\N}h_{11}^{(n)}\cdot (\pfdimtwo\obstuff{X})^{-n},
\end{gather*}
by assumption. Taking both together, the claim follows.
\end{proof}

The following is then key:

\begin{Lemma}\label{L:RTReductiveTwo}
Let $G=\mathrm{SL}_{2}(\C)$ and $\obstuff{X}$ any simple $G$-representation
of dimension $\dim_{\C}\obstuff{X}>1$. Then the growth problem 
$(G,\obstuff{X})$ is transient. The same is true when replacing $\C$ by any 
characteristic zero field.
\end{Lemma}

\begin{proof}
Recall that the fusion rules for $G=\mathrm{SL}_{2}(\C)$ are as follows. 
We index the simple $G$-representations $L(\lambda)$ 
by their highest weight $\lambda\in\N$. Then, 
using the convention that $L(\nu)=0$ if $\nu\notin\N$, we have
\begin{gather*}
L(\lambda)\otimes L(\mu)
\cong
L(\lambda-\mu)\oplus L(\lambda-\mu+2)\oplus\dots
\oplus L(\lambda+\mu-2)\oplus L(\lambda+\mu)
\end{gather*}
for $\lambda\geq\mu$ and similarly for $\lambda<\mu$. 

The fusion graphs are now easily computed. We give a few examples of how their 
action matrices look like:
\begin{gather}
M(\lambda=1)=
\begin{psmallmatrix}
0 & 1 & 0 & 0 & 0 & \dots
\\
1 & 0 & 1 & 0 & 0 & \dots
\\
0 & 1 & 0 & 1 & 0 & \dots
\\
0 & 0 & 1 & 0 & 1 & \dots
\\
0 & 0 & 0 & 1 & 0 & \dots
\\
\vdots & \vdots & \vdots & \vdots & \vdots & \vdots
\end{psmallmatrix}
,\quad
M(\lambda=3)=
\begin{psmallmatrix}
0 & 1 & 0 & 1 & 0 & \dots
\\
1 & 0 & 1 & 0 & 1 & \dots
\\
0 & 1 & 0 & 1 & 0 & \dots
\\
1 & 0 & 1 & 0 & 1 & \dots
\\
0 & 1 & 0 & 1 & 0 & \dots
\\
\vdots & \vdots & \vdots & \vdots & \vdots & \vdots
\end{psmallmatrix}
,
\end{gather}
and the number of offdiagonals with $1$s keep on 
increasing with steps of two for $\lambda=5$, $\lambda=7$ etc.
The case where $\lambda$ is even is similar, but with $1$s on the 
main diagonal except in the top left entry which is still zero. For example:
\begin{gather}
M(\lambda=2)=
\begin{psmallmatrix}
0 & 1 & 0 & 0 & 0 & \dots
\\
1 & 1 & 1 & 0 & 0 & \dots
\\
0 & 1 & 1 & 1 & 0 & \dots
\\
0 & 0 & 1 & 1 & 1 & \dots
\\
0 & 0 & 0 & 1 & 1 & \dots
\\
\vdots & \vdots & \vdots & \vdots & \vdots & \vdots
\end{psmallmatrix}
,\quad
M(\lambda=4)=
\begin{psmallmatrix}
0 & 1 & 1 & 0 & 0 & \dots
\\
1 & 1 & 1 & 1 & 0 & \dots
\\
1 & 1 & 1 & 1 & 1 & \dots
\\
0 & 1 & 1 & 1 & 1 & \dots
\\
0 & 0 & 1 & 1 & 1 & \dots
\\
\vdots & \vdots & \vdots & \vdots & \vdots & \vdots
\end{psmallmatrix}
.
\end{gather}
Note that the matrices for $\lambda\in 2\N$ are on a different basis 
on the Grothendieck level (namely on the Grothendieck classes of $L(0)$, $L(2)$, $L(4)$, etc.)

Let us now prove the case where $\lambda=2$ as an example. 
In this case $m_{11}^{n}$ is given by A005043 on \cite{Oeis}. In particular, 
$m_{11}^{(n)}\sim\frac{3^{3/2}}{8\sqrt{\pi}}\cdot n^{-3/2}\cdot 3^{n}$ and this implies
that
\begin{gather*}
\sum_{n\in\N}m_{11}^{(n)}\cdot 3^{-n}
\approx\tfrac{3^{3/2}}{8\sqrt{\pi}}
\sum_{n\in\N}n^{-3/2}<\infty,
\end{gather*}
which is the usual calculation one need to do in order to verify that a problem is 
transient. (For the reader unfamiliar with random walks, this is one of many possible, and equivalent, definitions of transient, see {\eg} \cite[Chapter 7]{Ki-symbolic-dynamics}.)

For general even $\lambda>0$ we have
$m_{11}^{(n)}\sim\text{some scalar}\cdot n^{-3/2}\cdot\lambda^{n}$, as follows 
from e.g. \cite[Theorem 2.2]{Bi-asymptotic-lie}. So the same calculations as for $\lambda=2$ works.
Moreover, again by e.g. \cite[Theorem 2.2]{Bi-asymptotic-lie}, when $\lambda$ is odd, then 
$m_{11}^{(2n)}\sim\text{some scalar}\cdot n^{-3/2}\cdot\lambda^{2n}$, and again the same 
calculation works.
\end{proof}

We can assume that $G$ is an algebraic group since we can replace it 
with the Zariski closure of its image in $\mathrm{GL}(\obstuff{X})$.
So assume that $G$ is an algebraic group, whose identity component is a reductive group since $\obstuff{X}$ is completely reducible. We can, and will, even replace $G$ by its identity component.

If $G$ is not a torus, 
then we can apply the Jacobson--Morozov theorem to find 
a copy of $\mathrm{SL}_{2}(\C)$ in $G$. Then \autoref{L:RTReductiveTwo} (see also \autoref{E:GrowthRateSLTwoThree})
and \autoref{L:RTReductiveOne}, whose assumptions are satisfied by \cite[Theorem 1.4]{CoOsTu-growth}, imply 
that the growth problem $(G,\obstuff{X})$ is transient.

On the other hand, if $G$ is torus, then Polya's classification of 
recurrent random walks applies, and we get a recurrent growth problem
if and only if the rank of the torus is $0$, $1$ or $2$.

The claim about positively recurrent versus null recurrent follows 
then from \autoref{E:GrowthRateSLTwoThree}.
\end{proof}

\begin{Proposition}\label{P:RTReductive}
The following growth problems are positively recurrent.
\begin{enumerate}

\item For an arbitrary field, $(G,\obstuff{X})$ for a finite group $G$ and $\obstuff{X}$ 
any $G$-representation.

\item For an arbitrary field, $(\catstuff{C},\obstuff{X})$ for a finite 
tensor category $\catstuff{C}$ and $\obstuff{X}\in\catstuff{C}$ any object.

\item For an arbitrary field, $(\sbim,\obstuff{X})$ for Soergel bimodules $\sbim$
attached to a finite Coxeter group and $\obstuff{X}\in\sbim$ any object.

\end{enumerate}
Moreover, the following growth problems are transient.
\begin{enumerate}[label=\emph{\upshape(\roman*)}]

\item For a field of characteristic zero, $(G,\obstuff{X})$ for a group $G$ 
and $\obstuff{X}$ any $G$-representation such that the
the Zariski closure of the image of 
$G$ in $\mathrm{GL}(\obstuff{X})$ is not virtually a torus.

\item In the below specified cases, $\big(U_{q},\obstuff{X}\big)$ for a quantum enveloping algebra 
$U_{q}=U_{q}(\mathfrak{g})$ (in the sense of \cite{Lu-qgroups-root-of-1} or \cite{AnPoWe-representation-qalgebras}) and $\obstuff{X}$ any nontrivial tilting $U_{q}(\mathfrak{g})$-representation
(meaning not a direct sum of one dimensional $U_{q}$-representations). The cases are:
\begin{enumerate}[label=$\bullet$]

\item $\K$ is arbitrary, and $q$ is not a root of unity.

\item $\K$ is an algebraically closed field of 
characteristic zero, and $q$ is a root of unity.
    
\end{enumerate}

\item For a field of characteristic zero, $(\sbim,\obstuff{X})$ for Soergel bimodules 
$\sbim$ attached to an affine Weyl group and $\obstuff{X}\in\sbim$ a generating object.

\end{enumerate}
\end{Proposition}

\begin{proof}
\text{(a).} A special case of (b).

\text{(b).} We can and will assume that $\obstuff{X}$ is a 
generating object: if that is not the case then we would go to a smaller tensor category. By \cite[Proposition A.1]{CoEtOsTu-sltwo-charp}, every projective indecomposable object of $\catstuff{C}$ is a direct summand of $\obstuff{X}^{\otimes d}$ for some $d\in\N$.

Let us first analyze the subgraph $\Gamma^{p}$ of the fusion 
graph $\Gamma$ that contains only vertices for the 
projective indecomposable objects. We call this the 
\emph{projective cell}, and we claim it 
is strongly connected and basic. To see this we first note that \cite[Proposition A.1]{CoEtOsTu-sltwo-charp} also gives us that the regular object appears in some tensor power 
of $\obstuff{X}$. Then \cite[Section 3.3]{EtGeNiOs-tensor-categories} implies that the 
projective cell is of maximal PF dimension since going from using Grothendieck classes of
simple objects to compute PF dimensions, as in \cite{EtGeNiOs-tensor-categories}, 
to indecomposable objects can only decrease the PF dimension. That the projective cell 
is strongly connected also follows from the existence of the regular objects in some tensor power of $\obstuff{X}$.

Moreover, recall that projectives form a $\otimes$-ideal, see e.g. \cite[Proposition 4.2.12]{EtGeNiOs-tensor-categories}. This also implies by \cite[Proposition A.1]{CoEtOsTu-sltwo-charp} that no strongly connected component without 
projective indecomposable objects is final since there is always a path 
to the projective cell.  

Since every finite tensor category 
has only finitely many indecomposable projective objects, 
this problem is positively recurrent.

\text{(c).} In this case $\Gamma$ is finite.

\text{(i).} As in the proof of \autoref{T:RTPolya}, 
the Jacobson--Morozov theorem implies that we will find
$\mathrm{SL}_{2}\subset\mathrm{GL}(\obstuff{X})$, and \autoref{L:RTReductiveOne}
and \autoref{L:RTReductiveTwo} then imply the claim.

\text{(ii).} If $q$ is not a root of unity, then this problem has the same combinatorics 
as (i), see for example \cite[Corollary 7.7]{AnPoWe-representation-qalgebras}, so we are done. In the other case there are projective $U_{q}$-representations 
and they appear in some tensor power of $\obstuff{X}$:

\begin{Lemma}\label{L:RTReductiveThree}
Every projective indecomposable tilting $U_{q}$-representation is a direct summand of $\obstuff{X}^{\otimes d}$ for some $d\in\N$.
\end{Lemma}

\begin{proof}
We use two results from \cite[Theorem 9.12]{AnPoWe-representation-qalgebras}: 
there are enough projective tilting objects, and every projective tilting object is also injective. With these two facts, the proof of the lemma is, mutatis mutandis, the 
same as \cite[Proof of Proposition A.1]{CoEtOsTu-sltwo-charp}.
\end{proof}

Using \autoref{L:RTReductiveThree}, we get the the only final class can be 
the projective cell, and it remains to argue that this cell is transient. The results 
in \cite[Theorem 3.1 and Remark 2.(2)]{An-steinberg-linkage} imply that the projective 
cell (called the Steinberg component in \cite{An-steinberg-linkage}) is a copy 
of the whole category in the semisimple case upon factoring the Steinberg $U_{q}$-representation. This implies 
that the projective cell is transient: The majority of paths in the starting category
will eventually be in the projective cell, so they move away from the origin. Then 
the equivalence in \cite[Theorem 3.1 and Remark 2.(2)]{An-steinberg-linkage} jumps back 
to the original problem for which we already know that the majority of paths eventually leave the origin.

\text{(iii).} Using the correspondence between the cells from 
\cite{Os-tensor-ideals-tilting}, 
very similar as for the quantum group at a complex root of unity. 
Details are omitted.
\end{proof}

\begin{Remark}\label{R:RTReductive}
One could guess that the analogs of \autoref{P:RTReductive}.(i) is also 
true in positive characteristic, e.g. for $G=\mathrm{SL}_{2}(\bF[p])$ this follows from \autoref{E:RecurrentSL2Fp}. Similarly, \autoref{P:RTReductive}.(ii) 
might also be true for quantum groups in the mixed case, and \autoref{P:RTReductive}.(iii) 
might be true for all infinite Coxeter types and regardless of the characteristic of the underlying field.
\end{Remark}


\section{Asymptotics for positively recurrent categories}\label{S:Recurrent}


Positively recurrent categories are essentially finite with respect 
to growth problems as justified by the main result of this section.
However, we need a special case of positively recurrent growth problems, see \autoref{D:RecurrentGrowth}.

\begin{Convention}
In this section, all growth problems that we consider are positively recurrent in the sense of 
\autoref{D:RTMain}.
\end{Convention}

First, let us consider the fusion graph $\Gamma$ 
associated to some growth problem $(R,C)$. Take the naive cutoff $(\Gamma_{k})_{k\in\N}$ and let $\lambda_{k}=\pfdimplain\Gamma_{k}$ (since $\Gamma_{k}$ is finite, this is the classical 
PF dimension.)
Let $\lambda^{sec}_{k}$ be any second largest eigenvalue of $\Gamma_{k}$. 

\begin{Lemma}\label{L:RecurrentGrowthExistence}
We have $\lambda=\lim_{k\to\infty}\lambda_{k}\in\Rplus$ and $\lambda^{sec}=\lim_{k\to\infty}\lambda^{sec}_{k}\in\C$.
\end{Lemma}

\begin{proof}
Note that $0\leq\lambda_{k}\leq\lambda_{k+1}$ so that it remains to argue why $\lambda=\lim_{k\to\infty}\lambda_{k}\neq\infty$.
As usual in the theory, see e.g. \cite[Section 7]{Ki-symbolic-dynamics}, $1/\lambda$ is the radius of 
convergence of a certain power series, which in turn has a positive radius of convergence on any FBC. 
Thus, we get $\lambda<\infty$. The second largest eigenvalue exists in $\Rplus$ by the same reasoning 
and $\lambda<\infty$ so $|\lambda^{sec}|<\infty$.
\end{proof}

Using \autoref{L:RecurrentGrowthExistence}, we will write $\lambda=\lim_{k\to\infty}\lambda_{k}$ and $\lambda^{sec}=\lim_{k\to\infty}\lambda^{sec}_{k}$ for a given 
positively recurrent growth problem.

\begin{Definition}\label{D:PFAnother}
We call $\lambda$ the (positively recurrent) \emph{PF dimension} of $(R,c)$, and write $\pfdimthree c=\lambda$. (Here we always use the naive cutoff.)
\end{Definition}

%

\begin{Definition}\label{D:RecurrentGrowth}
Let $(R,c)$ be a positively recurrent 
growth problem as defined in \autoref{D:RTMain}. 
We say 
$(R,c)$ is \emph{sustainably positively recurrent} if:
\begin{enumerate}

\item There is only one basic class.

\item There exists $\varepsilon>0$ such that
the PF dimensions of the nonbasic classes are smaller than $\lambda-\varepsilon$.

\item We have $(n\mapsto m_{ij}^{(n)})\in o\big((\lambda+\varepsilon)^{n}\big)$ for all $\varepsilon\in\R_{>0}$.

\end{enumerate}
We also say positively recurrent and sustainable instead of 
sustainably positively recurrent.
\end{Definition}

We will denote the basic class from \autoref{D:RecurrentGrowth}.(a) 
by $C^{FBC}(\Gamma)$, which necessarily positive recurrent and a final basic class. However, $C^{FBC}(\Gamma)$ is not assumed to be finite or final in $\Gamma$.

\begin{Example}\label{E:RecurrentGrowthFinite}
All finite growth problems are sustainably positively recurrent. 
This includes all positively recurrent growth problems in \autoref{T:RTPolya}.
\end{Example}

Here are more examples of positively recurrent and sustainable growth problems:

\begin{Proposition}\label{P:RTReductiveSustainable}
All the positively recurrent examples in \autoref{P:RTReductive} are sustainable. In all these cases we have $\pfdimthree\obstuff{X}\in\R_{\geq 1}$, and $C^{FBC}(\Gamma)$ is finite and final in $\Gamma$.
\end{Proposition}

\begin{proof}
We have $\pfdimthree\obstuff{X}\in\R_{\geq 1}$ since the fusion graphs have edge weights in $\Z_{\geq 1}$. This follows from the argument in \autoref{E:GrowthRatePFAdmissible} together with the observation that every vertex admits a path to the projective cell.

Second, the growth problem in (a) is a special case of the one in (b), while 
the growth problem in (c) has a finite fusion graph. So it remains to go through in 
\autoref{D:RecurrentGrowth} (a)-(c) for case (b), which we do after the following lemma.

\begin{Lemma}\label{L:RTTensorCategoryGrowth}
For any growth problem as in 
\autoref{P:RTReductive}.(b) we have 
\begin{gather*}
\bsymbol(n)\sim f(n)\cdot(\pfdimthree\obstuff{X})^{n}\text{ with }f\colon\N\to(0,1]\text{ periodic with finite period}.
\end{gather*}
\end{Lemma}

\begin{proof}
In this case let $\lsymbol_{n}=\ell(\obstuff{X}^{\otimes n})$ where 
$\ell$ is the length. The fusion graph of the growth problem associated to the length 
is finite and it follows (e.g. by copying \cite{LaTuVa-growth-pfdim}) that
\begin{gather*}
\lsymbol(n)\sim\tilde{f}(n)\cdot(\pfdimplain\obstuff{X})^{n}\text{ with }\tilde{f}\colon\N\to(0,1]\text{ periodic with finite period},
\end{gather*}
where $\pfdimplain\obstuff{X}$ is the usual PF dimension of this finite
growth problem. Note that
\begin{gather*}
\bsymbol_{n}\leq\lsymbol_{n}
\end{gather*}
since $\bsymbol_{n}$ is defined to count indecomposable summands, while $\lsymbol_{n}$ counts simple objects
in the Jordan--H{\"o}lder filtration. Then \cite[Proposition A.1]{CoEtOsTu-sltwo-charp} implies that 
this $\pfdimthree\obstuff{X}$ is the one from \autoref{L:RecurrentGrowthExistence}.

As the next step, using the 
universal grading group of $\catstuff{C}$, one can see that
$b(n)/(\pfdimthree\obstuff{X})^{n}$ has only finitely many limiting points. This implies that the generating function of $b_{n}$ satisfies the properties necessary to run e.g. 
\cite[Section 7.7]{Mi-analytic-combinatorics}, 
and we get $b(hn+s)\sim t_{h,s}\cdot(\pfdimthree\obstuff{X})^{hn}$ for some $t_{h,s}\in(0,1]$, where $h$ is the number of limiting points and $s\in\{0,\dots,h-1\}$. Collecting the scalars $t_{h,s}\in(0,1]$ into a piecewise constant function $f\colon\N\to[0,1)$ shows the claim.
\end{proof}

\textit{Property (a).} As in the proof of 
\autoref{P:RTReductive}, by \cite[Proposition A.1]{CoEtOsTu-sltwo-charp} the unique FBC is the projective cell, and its 
PF dimension is $\pfdimthree\obstuff{X}$.
We additionally argue now that all other strongly connected
components have PF dimension $<\pfdimthree\obstuff{X}$. Indeed, assuming the contrary leads to a contradiction with \autoref{L:RTTensorCategoryGrowth}. This is easy to see, for example in the toy case
\begin{gather*}
\Gamma=
\begin{tikzpicture}[anchorbase]
\begin{scope}[every node/.style={circle,inner sep=0pt,text width=6mm,align=center,draw=black,fill=white}]
\node (A) at (0,0) {};
\node (B) at (2,0) {};
\end{scope}
\begin{scope}[>={Stealth[black]},
every edge/.style={draw=orchid,very thick}]
\path [->] (A) edge (B);
\path [->] (A) edge [loop left] node{$\pfdimthree\obstuff{X}$} (A);
\path [->] (B) edge [loop right] node{$\pfdimthree\obstuff{X}$} (B);
\end{scope}
\end{tikzpicture}
,\quad
M(\Gamma)
=\begin{psmallmatrix}
\pfdimthree\obstuff{X} & 0
\\
1 & \pfdimthree\obstuff{X}
\end{psmallmatrix}
,
\end{gather*}
one already gets $b(n)\sim\frac{n}{\pfdimthree\obstuff{X}}\cdot(\pfdimthree\obstuff{X})^{n}$.
Recalling that all strongly connected components have some path to the projective cell,
the general case then reduces to at least this toy case.

\textit{Property (b).} Similarly as in (a), having 
basic classes $B_{n}$ with 
$\lim_{n\to\infty}\pfdimtwo B_{n}=\pfdimthree\obstuff{X}$ contradicts \autoref{L:RTTensorCategoryGrowth}. For example,
\begin{gather*}
\Gamma(Kl^{\prime})=
\begin{tikzpicture}[anchorbase]
\begin{scope}[every node/.style={circle,inner sep=0pt,text width=6mm,align=center,draw=black,fill=white}]
\node (A) at (0,0) {};
\node (B) at (2,0) {};
\node (C) at (4,0) {};
\node (D) at (6,0) {};
\node (E) at (8,0) {};
\node (F) at (10,0) {};
\node (P) at (6,-2) {};
\end{scope}
\node (G) at (12,0) {$\dots$};
\begin{scope}[>={Stealth[black]},
every edge/.style={draw=orchid,very thick}]
\path [->] (A) edge (B);
\path [->] (B) edge (C);
\path [->] (C) edge (D);
\path [->] (D) edge (E);
\path [->] (E) edge (F);
\path [->] (F) edge (G); 
\path [->] (B) edge (P);
\path [->] (C) edge (P);
\path [->] (D) edge (P);
\path [->] (E) edge (P);
\path [->] (F) edge (P);
\path [->] (G) edge (P);
\path [->] (P) edge [loop below] node {$1$} (P);
\path [->] (A) edge [loop above] node{$1-1$} (A);
\path [->] (B) edge [loop above] node{$1-1/2$} (B);
\path [->] (C) edge [loop above] node{$1-1/3$} (C);
\path [->] (D) edge [loop above] node{$1-1/4$} (D);
\path [->] (E) edge [loop above] node{$1-1/5$} (E);
\path [->] (F) edge [loop above] node{$1-1/6$} (F);
\end{scope}
\end{tikzpicture}
,
\end{gather*}
has $b(n)\in\Omega(1.5^{n})$ while $\pfdimthree\obstuff{X}=1$.

\textit{Property (c).}
As in the points above, nonsustainable would contradict \autoref{L:RTTensorCategoryGrowth} so we are done.
Roughly, consider the following modification of the graph from \autoref{Eq:GrowthRateMainTheoremKlein}:
\begin{gather*}
\Gamma(Kl^{\prime\prime})=
\begin{tikzpicture}[anchorbase]
\begin{scope}[every node/.style={circle,inner sep=0pt,text width=6mm,align=center,draw=black,fill=white}]
\node (A) at (0,0) {};
\node (B) at (2,0) {};
\node (C) at (4,0) {};
\node (D) at (6,0) {};
\node (E) at (8,0) {};
\node (F) at (10,0) {};
\node (P) at (6,-2) {};
\end{scope}
\node (G) at (12,0) {$\dots$};
\begin{scope}[>={Stealth[black]},
every edge/.style={draw=orchid,very thick}]
\path [->] (A) edge (B);
\path [->] (B) edge (C);
\path [->] (C) edge (D);
\path [->] (D) edge (E);
\path [->] (E) edge (F);
\path [->] (F) edge (G); 
\path [->] (B) edge node[left,xshift=-0.16cm]{$2^{1}$} (P);
\path [->] (C) edge node[left,xshift=-0.15cm]{$2^{2}$} (P);
\path [->] (D) edge node[left,xshift=-0.15cm]{$2^{3}$} (P);
\path [->] (E) edge node[left,xshift=-0.15cm]{$2^{4}$} (P);
\path [->] (F) edge node[left,xshift=-0.15cm]{$2^{5}$} (P);
\path [->] (G) edge node[left,xshift=-0.15cm]{$\dots$} (P);
\path [->] (P) edge [loop below] node{$1$} (P);
\end{scope}
\end{tikzpicture}
\end{gather*}
The growth of $b_{n}$ in this case overshoots $(\pfdimthree\obstuff{X})^{n}=1^{n}$ and we have 
$b(n)\in\Omega(1.5^{n})$.

\textit{Final and finite in $\Gamma$.} We already discussed this above.
\end{proof}

Fix the naive cutoff $\Gamma_{k}$. Let $h_{k}$ be the period and $\zeta_{k}=\exp(2\pi i/h_{k})$.
Similarly as in \cite[Section 1]{LaTuVa-growth-pfdim},
let us denote the right (the one for the left action) and left
(the one for the right action) eigenvectors for $\zeta_{k}^{i}\lambda_{k}$ by 
$v_{i}^{k}$ and $w_{i}^{k}$, normalized such that $(w^{k}_{i})^{T}v_{i}^{k}=1$.
Let $v_{i}^{k}(w_{i}^{k})^{T}[1]$ denote taking the sum of the first column of the matrix 
$v_{i}^{k}(w_{i}^{k})^{T}$. Using this, we define 
\begin{gather*}
\asymbol_{k}(n)=
\big(v_{0}^{k}(w_{0}^{k})^{T}[1]\cdot 1+v_{1}^{k}(w_{1}^{k})^{T}[1]\cdot\zeta^{n}+v_{2}^{k}(w_{2}^{k})^{T}[1]\cdot(\zeta^{2})^{n}+\dots+v_{h-1}^{k}(w_{h-1}^{k})^{T}[1]\cdot(\zeta^{h-1})^{n}\big)
\cdot\lambda^{n}.
\end{gather*}
Similarly, but directly for $\Gamma$, we can define 
\begin{gather}\label{Eq:RecurrentMainFormula}
\asymbol(n)=
\big(v_{0}(w_{0})^{T}[1]\cdot 1+v_{1}(w_{1})^{T}[1]\cdot\zeta^{n}+v_{2}(w_{2})^{T}[1]\cdot(\zeta^{2})^{n}+\dots+v_{h-1}(w_{h-1})^{T}[1]\cdot(\zeta^{h-1})^{n}\big)
\cdot\lambda^{n},
\end{gather}
where $h$ is the period of the FBC of $\Gamma$ (for strongly connected graphs like the FBC this is defined in e.g. \cite[Begin of Section 7.1]{Ki-symbolic-dynamics}). Note that,
a priori, $a(n)$ might not be a finite expression.

\begin{Theorem}\label{T:RecurrentGrowth}
Assume that the growth problem $(R,C)$ is sustainably positively recurrent. Then:
\begin{enumerate}

\item We have $a(n)\in\Rplus$ and
\begin{gather*}
\bsymbol(n)\sim\asymbol(n).
\end{gather*}

\item If the FBC is final in $\Gamma$, 
the limit $\lim_{k\to\infty}\asymbol_{k}(n)$ exists and is equal to $\asymbol(n)$.

\item If the FBC is finite and final in $\Gamma$, 
then the convergence of $\lim_{n\to\infty}b_{n}/a_{n}=1$ is geometric with ratio $|\lambda^{sec}/\lambda|$ and $|b_{n}-a_{n}|\in O\big((\lambda^{sec})^{n}+n^{d}\big)$ for some $d\in\R_{>0}$.

\end{enumerate}
\end{Theorem}

\begin{Remark}\label{R:RecurrentGrowthLazy}
\cite{LaTuVa-growth-pfdim} states an analog of \autoref{T:RecurrentGrowth} for 
arbitrary finite graphs, but the statement itself is a bit nasty. We decided not to include its positively recurrent 
analog in this paper.
\end{Remark}

\begin{proof}[Proof of \autoref{T:RecurrentGrowth}]
We start with the statement about the asymptotic and then proof the other claims that we call ``Finite approximation'' and ``Variance''.

\textit{Asymptotic.}
Condition (a) of \autoref{D:RecurrentGrowth} implies that there is a unique 
basic class $C^{FBC}(\Gamma)$ (which is a final basic class, of course). Note that 
$\lambda<\infty$ by \autoref{L:RecurrentGrowthExistence}.
Moreover, $\lambda=0$ if and only if the growth problem is 
zero, and hence we can and will assume that $\lambda\in\R_{>0}$.

We borrow from and adjust \cite[Section 7.1]{Ki-symbolic-dynamics}. There are two differences between our setting and \cite[Section 7.1]{Ki-symbolic-dynamics} to keep in mind: Firstly, and crucial, we consider graphs that are potentially not strongly connected. Secondly, we can have a period, but that is a rather harmless generalization. We postpone this discussion to the end of the proof and assume first that we are in the aperiodic situation.

Unless stated otherwise, our indexing set below will be the vertices of the fusion graph $\Gamma$ associated to $(R,C)$, and these vertices are often called $i,j$.

\begin{Notation}\label{N:RecurrentGrowth}
To have the same conventions as \cite{Ki-symbolic-dynamics} in this proof, let $T$ be the transpose of $M$.
\end{Notation}

We define three generating functions $T$, $L$ and $R$:
\begin{gather*}
T_{ij}(z)=\sum_{n\in\N}t_{ij}^{(n)}z^{n}
,\quad
L_{ij}(z)=\sum_{n\in\N}l_{ij}(n)z^{n}
,\quad
R_{ij}(z)=\sum_{n\in\N}r_{ij}(n)z^{n}
,
\\
L_{ij}^{\prime}=\sum_{n\in\N}nl_{ij}(n)z^{n}
,\quad
\mu(i)=\lambda^{-1}L_{ii}^{\prime}(\lambda^{-1})\text{ (with $\mu(i)=\infty$ allowed)},
\end{gather*}
where the coefficients of $L$ and $R$ are defined as follows.
The coefficient $l_{ij}(n)$ is the sum of the labels of the paths 
$i\to j$ of length $n$ that do not return to $i$ with $<n$ steps.
Similarly, $r_{ij}(n)$ is the sum of the labels of the paths 
$i\to j$ of length $n$ that do not return to $j$ with $<n$ steps.

\textit{Step 1.} We start with two crucial lemmas.
We define $\lambda_{ij}=
\limsup_{n\to\infty}\sqrt[n]{t_{ij}^{(n)}}$.

\begin{Lemma}\label{L:RecurrentInequalityPFDim}
We have $\lambda_{ij}\leq\lambda$ for all $i,j\in\Gamma$.
\end{Lemma}

\begin{proof}
Immediate from $(n\mapsto t_{ij}^{(n)})\in o\big((\lambda+\varepsilon)^{n}\big)$ for all $\varepsilon\in\R_{>0}$.
\end{proof} 

\begin{Lemma}\label{L:RecurrentInequalityPFDimTwo}
For all $i,j\in\Gamma$, the radii of convergences of 
$T_{ij}(z)$, $L_{ij}(z)$ and $R_{ij}(z)$ are $\geq\lambda^{-1}$.
\end{Lemma}

\begin{proof}
The radius of convergence of $T_{ij}(z)$ is $\lambda_{ij}^{-1}$ by the Cauchy--Hadamard theorem, and the result follows from \autoref{L:RecurrentInequalityPFDim} and the inequalities $l_{ij}(n)\leq t_{ij}^{(n)}$ and $r_{ij}(n)\leq t_{ij}^{(n)}$.
\end{proof}

We are now ready to go through statements in \cite[Section 7.1]{Ki-symbolic-dynamics}. The references below are all with respect to
\cite{Ki-symbolic-dynamics}.

\textit{Step 2.} Most parts of Lemma 7.1.6 are formal and
work verbatim. Some of the arguments need 
\autoref{L:RecurrentInequalityPFDimTwo}, but then can be proven 
mutatis mutandis with it.

\textit{Step 3.} Next, we see a mild modification. Namely, Lemma 7.1.7 for $i,j\in\Gamma$ in the same component works mutatis mutandis since 
$\lambda$ is the leading eigenvalue so $\lambda^{-1}$ is the minimal radius of convergence, see \autoref{L:RecurrentInequalityPFDimTwo}.
Moreover, $L_{ij}(\lambda^{-1})$ is finite in general if there is a path from $j$ to $i$, and vice versa for $R_{ij}(\lambda^{-1})$, by the same arguments as in \cite[Section 7.1]{Ki-symbolic-dynamics}.

\textit{Step 4.} Lemma 7.1.8 remains true with the following changes. If $i\in C^{FBC}(\Gamma)$, then $T_{ij}(\lambda^{-1})=\infty$ and $L_{ii}(\lambda^{-1})=1$ since the FBC is recurrent. We do not need the case where $i\notin C^{FBC}(\Gamma)$.

\textit{Step 5.} Choose some enumeration of the vertices of $\Gamma$. For each $i\in C^{FBC}(\Gamma)$ define a row and a column vector
\begin{gather*}
\ell^{(i)}=
\begin{pmatrix}
L_{i1}(\lambda^{-1}) & L_{i2}(\lambda^{-1}) & \dots & L_{ij}(\lambda^{-1}) & \dots &
\end{pmatrix},
\\
r^{(i)}=
\begin{pmatrix}
R_{1i}(\lambda^{-1}) & R_{2i}(\lambda^{-1}) & \dots & R_{ji}(\lambda^{-1}) & \dots & 
\end{pmatrix}^{T}.
\end{gather*}
These are well-defined by Step 3.

\textit{Step 6.} Lemma 7.1.9.(i) remains true, that is, the 
$\ell^{(i)}$ and the $r^{(i)}$ are left and right $\lambda$-eigenvectors 
of the action matrix $T$:
\begin{gather*}
\ell^{(i)}T=\lambda\cdot\ell^{(i)},
\quad
Tr^{(i)}=\lambda\cdot r^{(i)}.
\end{gather*}
(For the meticulous reader who want to double check that $r^{(i)}$ is an eigenvector we point out that the 
$R_{ii}$ in Lemma 7.1.6.(viii) should be a $R_{jj}$.)

Choose any $i\in C^{FBC}(\Gamma)$ and define
\begin{gather*}
\ell=\ell^{(i)},\quad r=r^{(i)}.
\end{gather*}
It is easy to see that the choice of the defining $i\in C^{FBC}(\Gamma)$ is not 
important for what we will do below.

\textit{Step 7.} Lemma 7.1.10 remains true, with the same proof, but needs some adjustments. Firstly, the vectors $x$ and $y$
are assumed to have a nonzero entry corresponding 
to some $k\in C^{FBC}(\Gamma)$. Then $x$ and $y$ are 
strictly positive on the entries corresponding 
to all $k\in C^{FBC}(\Gamma)$. This uses irreducibility of the 
matrix supported on $C^{FBC}(\Gamma)$. (ii) and (iii) then remain true without 
further change.

\textit{Step 8.}
For Lemma 7.1.11 we replace $xT\leq\lambda x$ and $Ty\leq\lambda y$ by $xT=\lambda x$ and $Ty=\lambda y$ 
and we keep the same assumptions on $x$ and $y$ from the previous step. We renormalize $x$ and $y$ such that $x_{i}=1$ and $y_{i}=1$, and the previous step ensure that $x-l^{(i)}$ and $y-r^{(i)}$ are nonnegative and have a zero entry on $C^{FBC}(\Gamma)$. We then obtain that $x-l^{(i)}$ and $y-r^{(i)}$ are supported on the nonbasic classes. We can extract from $x-l^{(i)}$ and $y-r^{(i)}$ eigenvectors of a matrix whose spectrum does not contain $\lambda$, because of the condition (b) in \autoref{D:RecurrentGrowth}. We finally obtain that $x=l^{(i)}$ and $y=r^{(i)}$.

\textit{Step 9.} Since the FBC is positively recurrent, we get 
$L_{ii}^{\prime}(\lambda^{-1})<\infty$ for $i\in C^{FBC}(\Gamma)$. It follows that $l\cdot r=\mu(i)<\infty$ (dot product) by copying the calculation in the proof of Lemma 7.1.14.

\textit{Step 10.} Lemma 7.1.15 remains true (it is independent of the growth problem).

\textit{Step 11.} Theorem 7.1.18 has two cases: for $i\notin C^{FBC}(\Gamma)$ we are in case (i), and for $i\in C^{FBC}(\Gamma)$
we are in case (ii). The justification for $i\notin C^{FBC}(\Gamma)$ is clear, and the proof for $i\in C^{FBC}(\Gamma)$ works mutatis mutandis as in \cite{Ki-symbolic-dynamics}, using the adjustments made in the steps above.

\textit{Step 12.} Lemma 7.1.19.(i) works verbatim. For Lemma 7.1.19.(ii) there are three cases. If $i\in C^{FBC}(\Gamma)$, then the first equality holds. Moreover, for $j\in C^{FBC}(\Gamma)$ the equality between the first and the last term holds. Thus, for $i,j\in C^{FBC}(\Gamma)$ Lemma 7.1.19.(ii) 
works in the same way as in \cite{Ki-symbolic-dynamics}.

\textit{Step 13.} Lemma 7.1.20 works with the extra assumption that $i,k\in C^{FBC}(\Gamma)$ in part (a) of (i), $j,k\in C^{FBC}(\Gamma)$ in part (b) of (i) and $i,j\in C^{FBC}(\Gamma)$ in (ii).

\textit{Step 14.} As in \cite{Ki-symbolic-dynamics}, steps 1-13 prove \cite[Theorem 7.1.3.(f)]{Ki-symbolic-dynamics} for sustainably positively recurrent growth problems.

\textit{Step 15.} Now, we can mimic the proof of \cite[Theorem 7]{LaTuVa-growth-pfdim} to obtain the precise asymptotic on $b(n)$.

For periodic matrices evaluate the power series at $\zeta^{k}\lambda$ for $\zeta$ a root of unity. With this change all the above work mutatis mutandis. See also the final part of 
\cite[Section 7.1]{Ki-symbolic-dynamics}.

\textit{Finite approximation.} Since the FBC is final, the vector $\ell$ is supported by the FBC. The first column of the matrices in \autoref{Eq:RecurrentMainFormula} are then supported by the FBC. 
Hence, \cite[Theorem 7.1.4]{Ki-symbolic-dynamics} applies and we are done.

\textit{Variance.} As in the previous point, but since the FBC is now assumed to be finite, we can use \cite[Theorem 1]{LaTuVa-growth-pfdim}.
\end{proof}

\begin{Remark}\label{R:RecurrentFinite}
For Krull--Schmidt monoidal categories with finitely many indecomposable objects
\autoref{T:RecurrentGrowth} recovers 
the main results of \cite{LaTuVa-growth-pfdim} 
(up to the point made in \autoref{R:RecurrentGrowthLazy}).
\end{Remark}

\begin{Example}\label{E:RecurrentKlein}
Let us come back to \autoref{E:GrowthRateFinite}.
The random walk on the fusion graph $\Gamma$ in 
\autoref{Eq:GrowthRateMainTheoremKlein} is positively recurrent.
Thus, \autoref{T:RecurrentGrowth} applies and we get
\begin{gather*}
\bsymbol_{n}\sim\tfrac{1}{4}\cdot 3^{n}
\end{gather*}
(note that $|G|=4$ and $\sum_{\text{$L$ simples}}\dim_{\F[2]}L=1$)
using the finite cutoffs as indicated in \autoref{Eq:GrowthRateMainTheoremKlein}.
The second largest eigenvalue for all finite cutoffs is zero, so the variance is bounded by some polynomial function. Indeed, we get:
\begin{gather*}
\begin{tikzpicture}[anchorbase]
\node at (0,0) {\includegraphics[height=4.6cm]{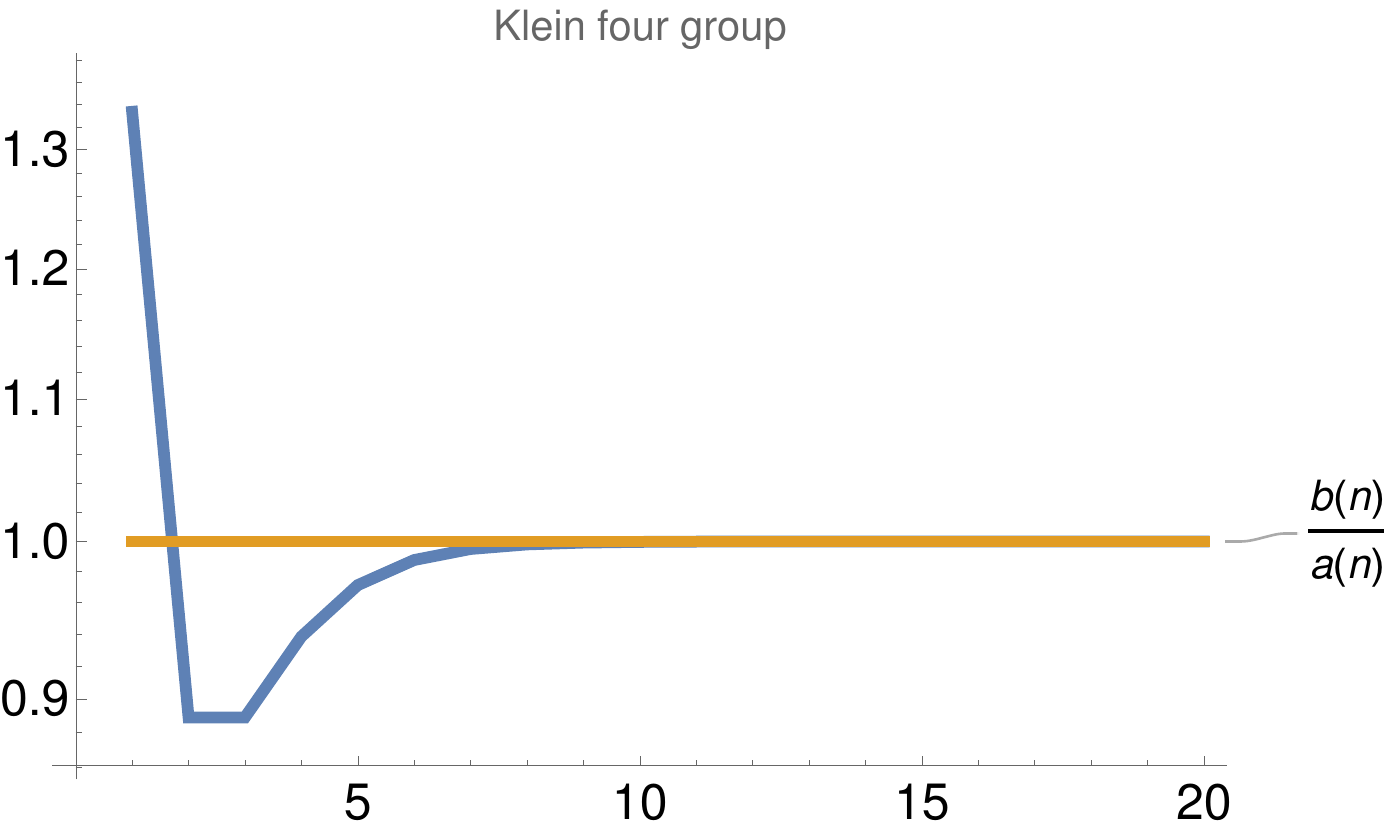}};
\end{tikzpicture}
,
\begin{tikzpicture}[anchorbase]
\node at (0,0) {\includegraphics[height=4.6cm]{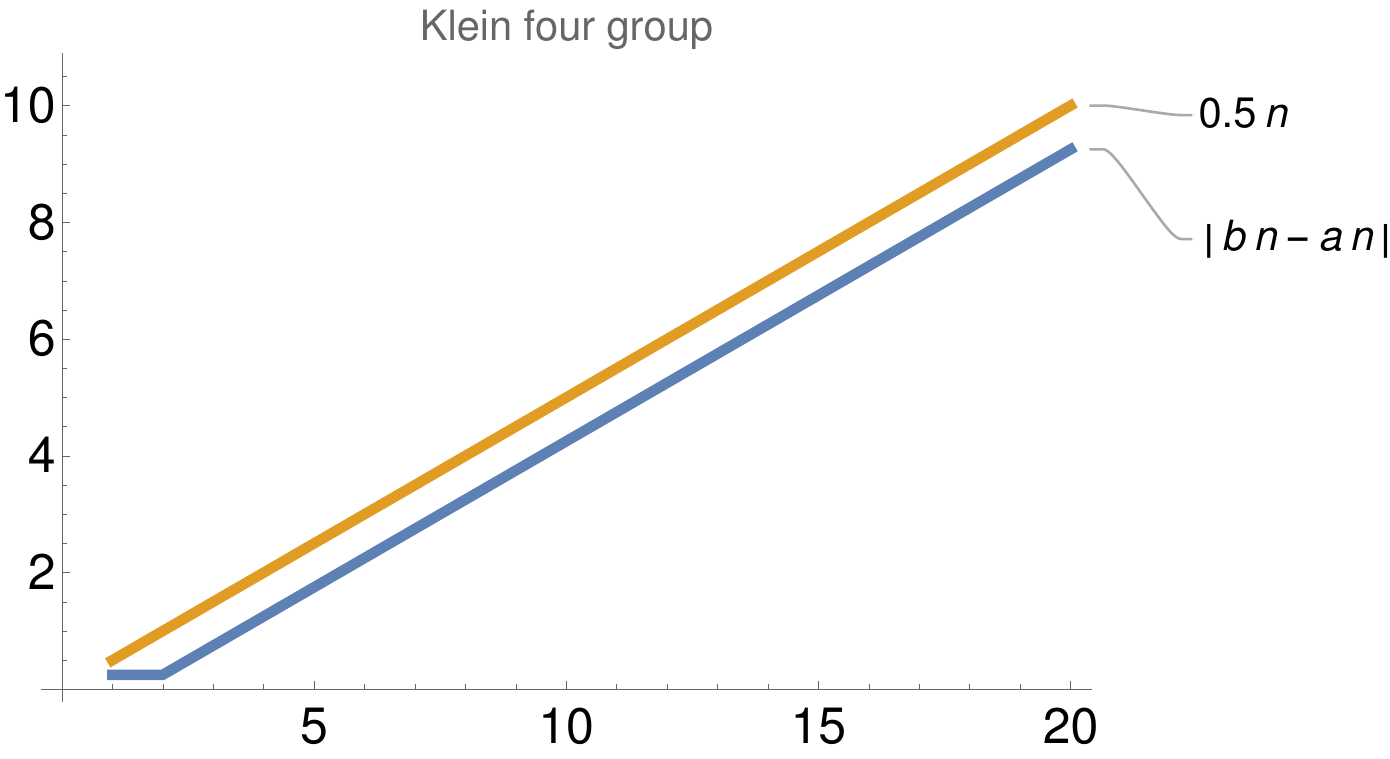}};
\end{tikzpicture}
\hspace{-0.5cm}.
\end{gather*}
The left plot is a logplot, the right is a standard plot.
\end{Example}

\begin{Example}\label{E:RecurrentPSL2}
Another example from the realm of finite groups is the following. Consider $G=\mathrm{PSL}_{2}(\F[7])$. Over $\F[2]$ the representation theory of $G$ is not 
semisimple. In this case there are two three dimensional simple representations, and let 
$\obstuff{X}$ be any of these two. The choice of this example is motivated
by \cite[Theorem 1.4]{Cr-tensor-simple-modules}, and is the smallest example on that list.

The fusion graph of the associated growth problem $(G,\obstuff{X})$ can be seen to be (illustrated below with a cutoff)
\begin{gather}\label{Eq:RecurrentPSL2}
\begin{gathered}
\begin{tikzpicture}[anchorbase]
\node at (0,0) {\includegraphics[height=4.5cm]{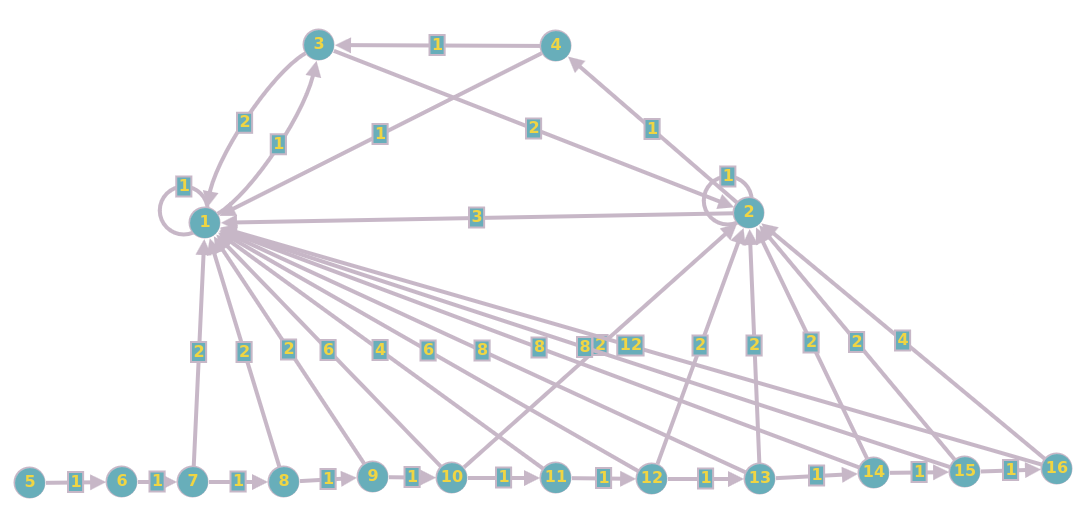}};
\end{tikzpicture}
,
\end{gathered}
\end{gather}
\begin{gather*}
\begin{gathered}
M(PSL2)=
\scalebox{0.5}{$\left(
\begin{array}{cccccccccccccccc}
\cellcolor{spinach!25}1 & \cellcolor{spinach!25}3 & \cellcolor{spinach!25}2 & \cellcolor{spinach!25}1 & 0 & 0 & 2 & 2 & 2 & 6 & 4 & 6 & 8 & 8 & 8 & 12 \\
\cellcolor{spinach!25}0 & \cellcolor{spinach!25}1 & \cellcolor{spinach!25}2 & \cellcolor{spinach!25}0 & 0 & 0 & 0 & 0 & 0 & 2 & 0 & 2 & 4 & 4 & 4 & 4 \\
\cellcolor{spinach!25}1 & \cellcolor{spinach!25}0 & \cellcolor{spinach!25}0 & \cellcolor{spinach!25}1 & 0 & 0 & 0 & 0 & 0 & 0 & 0 & 0 & 0 & 0 & 0 & 0 \\
\cellcolor{spinach!25}0 & \cellcolor{spinach!25}1 & \cellcolor{spinach!25}0 & \cellcolor{spinach!25}0 & 0 & 0 & 0 & 0 & 0 & 0 & 0 & 0 & 0 & 0 & 0 & 0 \\
0 & 0 & 0 & 0 & 0 & 0 & 0 & 0 & 0 & 0 & 0 & 0 & 0 & 0 & 0 & 0 \\
0 & 0 & 0 & 0 & 1 & 0 & 0 & 0 & 0 & 0 & 0 & 0 & 0 & 0 & 0 & 0 \\
0 & 0 & 0 & 0 & 0 & 1 & 0 & 0 & 0 & 0 & 0 & 0 & 0 & 0 & 0 & 0 \\
0 & 0 & 0 & 0 & 0 & 0 & 1 & 0 & 0 & 0 & 0 & 0 & 0 & 0 & 0 & 0 \\
0 & 0 & 0 & 0 & 0 & 0 & 0 & 1 & 0 & 0 & 0 & 0 & 0 & 0 & 0 & 0 \\
0 & 0 & 0 & 0 & 0 & 0 & 0 & 0 & 1 & 0 & 0 & 0 & 0 & 0 & 0 & 0 \\
0 & 0 & 0 & 0 & 0 & 0 & 0 & 0 & 0 & 1 & 0 & 0 & 0 & 0 & 0 & 0 \\
0 & 0 & 0 & 0 & 0 & 0 & 0 & 0 & 0 & 0 & 0 & 0 & 0 & 0 & 0 & 0 \\
0 & 0 & 0 & 0 & 0 & 0 & 0 & 0 & 0 & 0 & 0 & 1 & 0 & 0 & 0 & 0 \\
0 & 0 & 0 & 0 & 0 & 0 & 0 & 0 & 0 & 0 & 0 & 0 & 1 & 0 & 0 & 0 \\
0 & 0 & 0 & 0 & 0 & 0 & 0 & 0 & 0 & 0 & 0 & 0 & 0 & 1 & 0 & 0 \\
0 & 0 & 0 & 0 & 0 & 0 & 0 & 0 & 0 & 0 & 0 & 0 & 0 & 0 & 1 & 0 \\
\end{array}
\right)$}
.
\end{gathered}
\end{gather*}
The projective cell is illustrated at the top and has the adjacency matrix
\begin{gather*}
M=
\begin{psmallmatrix}
1 & 0 & 1 & 0 \\
3 & 1 & 0 & 1 \\
2 & 2 & 0 & 0 \\
1 & 0 & 1 & 0 \\
\end{psmallmatrix}
.
\end{gather*}
The dimensions of the projective indecomposables, in the order as in 
\autoref{Eq:RecurrentPSL2}, are $8$, $16$, $16$, $8$.
The eigenvalues of $M$ are $\{3,\frac{1}{2}(-1+i\sqrt{7}),\frac{1}{2}(-1-i\sqrt{7}),0\}$, and our dominating growth rate is $3$.

The dimensions of the indecomposables along the bottom in \autoref{Eq:RecurrentPSL2} is
\begin{gather*}
1,3,9,11,17,35,25,43,49,51,57,75,65,83,89,91,97,115,105,123,129,\dots,
\\
\text{pattern}\colon +6,+2,+6,+18,-10,+18,\text{repeat},
\\
\text{generating function}\colon
\frac{16x^{5}+5x^{4}+7x^{3}+8x^{2}+3x+1}{(x-1)^{2}(x^{3}+2x^{2}+2x+1)},
\end{gather*}
and this pattern continues. From this pattern one gets the whole graph.

Moreover, one can check that $v_{0}^{k}(w_{0}^{k})^{T}[1]$ for the leading eigenvectors 
converges to $\frac{15}{168}\approx 0.089$ as $k\to\infty$. 
(Note that $|G|=168$ and $\sum_{\text{$L$ simples}}\dim_{\F[2]}L=15$.)
Thus, \autoref{T:RecurrentGrowth} implies that:
\begin{gather*}
\bsymbol_{n}\sim\frac{15}{168}\cdot 3^{n}.
\\
\begin{tikzpicture}[anchorbase]
\node at (0,0) {\includegraphics[height=4.6cm]{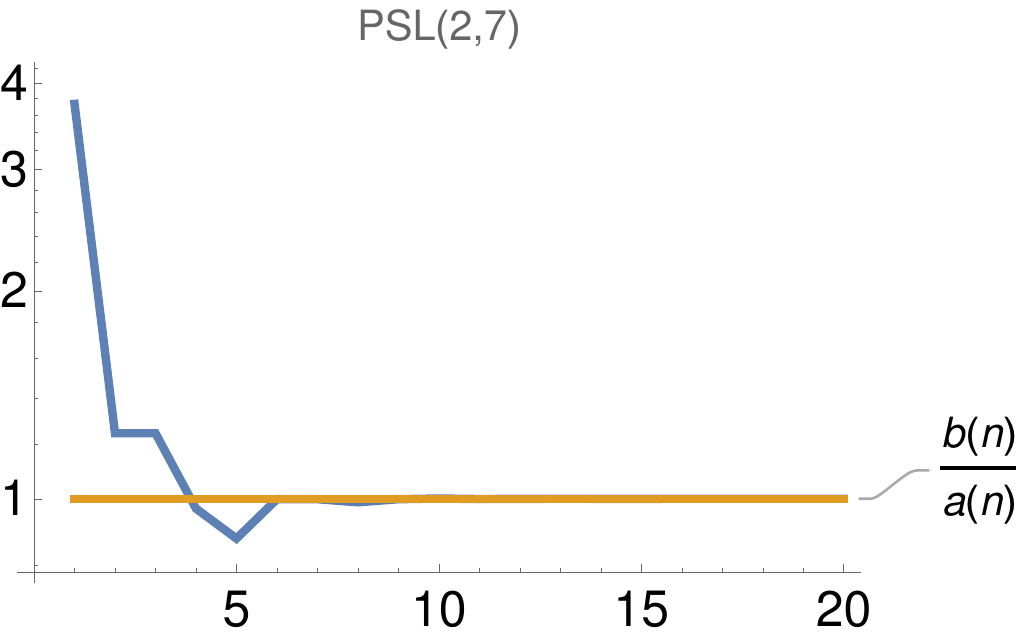}};
\end{tikzpicture}
,
\begin{tikzpicture}[anchorbase]
\node at (0,0) {\includegraphics[height=4.6cm]{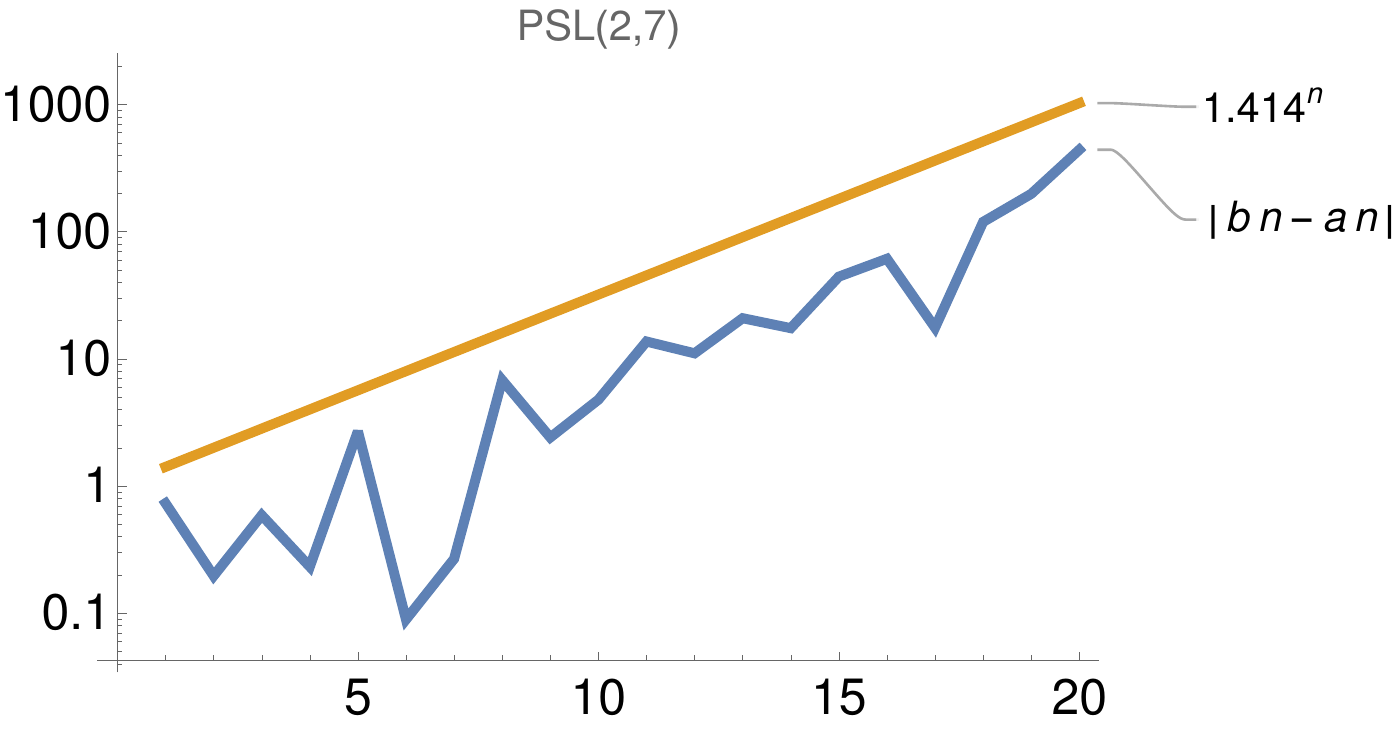}};
\end{tikzpicture}
\hspace{-0.5cm}.
\end{gather*}
Note that $|\lambda^{sec}|=|\frac{1}{2}(-1+i\sqrt{7})|=\sqrt{2}\approx 1.414$.
The convergence is geometric with ratio $\sqrt{2}/3=|\frac{1}{2}(-1+i\sqrt{7})|/3\approx 0.4714$.
\end{Example}

\begin{Example}\label{E:RecurrentSBim}
Let $W$ be a finite Coxeter group and $\catstuff{SBim}=\catstuff{SBim}(W)$ the associated 
monoidal category of Soergel bimodules. Let $\obstuff{X}\in\catstuff{SBim}$ 
be some generating object. The growth problem $(\catstuff{SBim},\obstuff{X})$
is positively recurrent and has the following asymptotics.

In the Grothendieck ring, the group ring $\Z W$, write
$[\obstuff{X}]=\sum_{s\in W}m_{s}\cdot s$ where $m_{s}\in\N$. One gets
\begin{gather*}
b(n)\sim\frac{1}{|W|}\cdot\big(\sum_{s\in W}m_{s}\big)^{n}.
\end{gather*}
This follows directly from \autoref{T:RecurrentGrowth} since one only has 
one FBC corresponding to the indecomposable Soergel bimodule associated with the longest word.
\end{Example}

\begin{Remark}
In \autoref{E:RecurrentSBim}, we lack a general expression for the second largest eigenvalue, except in the dihedral case for $\obstuff{X}$ being the indecomposable object corresponding to the product of the two simple reflections. For this growth problem, the eigenvalue depends on a parity condition related to whether the order of $W$ (always even) is divisible by four. Using, for example, \cite[Section 3]{Tu-sandwich-cellular} one can show that the coefficients of the minimal polynomial of the second largest eigenvalue in these cases are the signed versions of the sequences \cite[A085478]{Oeis} and \cite[A030528]{Oeis}, respectively. Thus, the second largest eigenvalues are closely related, though not identical, to the PF roots of the \emph{Morgan--Voyce polynomials}. A more general statement would be desirable.
\end{Remark}

\end{document}